\documentclass[final,leqno]{siamltex}
\usepackage[latin1]{inputenc}
\usepackage[T1]{fontenc}
\usepackage{times}
\usepackage[pdftex]{graphicx} 
\usepackage{tikz}

\expandafter\let\csname equation*\endcsname\relax
\expandafter\let\csname endequation*\endcsname\relax
\usepackage{amsmath}
\usepackage[english]{babel} 
\usepackage{graphicx}
\usepackage{cite}
\usepackage{amsfonts,amssymb}
\usepackage{color}  
\usepackage{epstopdf}
\usepackage{todonotes}
\usepackage{multirow}
\usepackage{verbatim}
\usepackage{hyperref}
\usepackage{enumerate}
\usepackage[ruled,vlined,linesnumbered]{algorithm2e}
\newtheorem{thm}{Theorem}
\newtheorem{lem}[theorem]{Lemma}
\newtheorem{prop}[theorem]{Proposition}
\newtheorem{prf}[theorem]{Proof}

\newcommand{\R}{{\mathbb R}}
\newcommand{\N}{{\mathbb N}}

\newcommand{\omegab}{\boldsymbol{\omega}}

\newcommand{\ts}[1]{\mbox{\textsf{#1}}}
\newcommand{\norm}[1]{\left\| #1 \right\|}

\newcommand{\ii}{{\mathrm{i}}}
\newcommand{\sigmab}{\boldsymbol{\sigma}}
\newcommand{\mub}{\boldsymbol{\mu}}

\newcommand{\xib}{\boldsymbol{\xi}}
\newcommand{\rhob}{\boldsymbol{\rho}}
\newcommand{\bm}[1]{\mathbf{#1}}
\newcommand{\C}{{\mathbb{C}}}
\newcommand{\Imag}{\mbox{Im}}
\newcommand{\Real}{\mbox{Re}}

\newcommand{\Sigmab}{\mathbf{\Sigma}}

\title{A variational non-linear constrained model for the inversion of FDEM data} 

\author{A Buccini$^1$ and P D\'iaz de Alba$^2$}

\begin{document}
\maketitle

	\begin{abstract}
		Reconstructing the structure of the soil using non-invasive techniques is a very relevant problem in many scientific fields, like geophysics and archaeology. This can be done, for instance, with the aid of Frequency Domain Electromagnetic (FDEM) induction devices. Inverting FDEM data is a very challenging inverse problem, as the problem is extremely ill-posed, i.e., sensible to the presence of noise in the measured data, and non-linear. Regularization methods substitute the original ill-posed problem with a well-posed one whose solution is an accurate approximation of the desired one. In this paper we develop a regularization method to invert FDEM data. We propose to determine the electrical conductivity of the ground by solving a variational problem. The minimized functional is made up by the sum of two term: the data fitting term ensures that the recovered solution fits the measured data, while the regularization term enforces sparsity on the Laplacian of the solution. The trade-off between the two terms is determined by the regularization parameter. This is achieved by minimizing an $\ell_2-\ell_q$ functional with $0<q\leq 2$. Since the functional we wish to minimize is non-convex, we show that the variational problem admits a solution. Moreover, we prove that, if the regularization parameter is tuned accordingly to the amount of noise present in the data, this model induces a regularization method. Some selected numerical examples on synthetic and real data show the good performances of our proposal.
	\end{abstract}

\section{Introduction}
In this paper we consider a severely ill-posed problem that arises in geophysics, namely the reconstruction of the electrical conductivity of the ground using a Frequency Domain Electromagnetic (FDEM) induction device. In this work we consider a two-dimensional vertical section of the ground and assume the magnetic permeability to be known and constant. The extension to the three-dimensional case and unknown magnetic permeability is straightforward and, for simplicity, we do not consider it here. The considered problem is of the form
\begin{equation}\label{eq:probl}
	\arg\min_{\Sigma\in\mathcal{X}}\norm{\mathcal{M}(\Sigma)-\mathcal{B}}_\mathcal{Y}^2,
\end{equation}
where $\mathcal{X}$ and $\mathcal{Y}$ are Hilbert spaces, $\Sigma\in\mathcal{X}$ represents the electrical conductivity, $\mathcal{B}\in\mathcal{Y}$ collects the (possibly noisy) data measured by the FDEM device, $\norm{\cdot}_\mathcal{Y}$ denotes the norm defined on $\mathcal{Y}$, and $\mathcal{M}:\mathcal{X}\rightarrow\mathcal{Y}$ is a non-linear function. We assume that \eqref{eq:probl} is ill-posed, i.e., the solution may not be unique and its computation is very sensible to the presence of noise in the measured data. Therefore, the naive solution of \eqref{eq:probl} is usually a poor approximation of the exact electrical conductivity. To compute a meaningful solution we need to resort to regularization methods; see, e.g., \cite{Engl} and references therein for a more detailed discussion. 

We represent the two-dimensional ground as the semi-infinite rectangle $[a,b]\times[0,\infty)$, where the first dimension is the horizontal space dimension (intuitively the ``sea-level'') and the second dimension represents the depth. Our purpose is to reconstruct an ``image'' of the electrical conductivity on this rectangle. This can be done using a FDEM induction device called Ground Conductivity Meter (GCM). Its principle of operation is based on an alternating electrical current which flows through a small electric wire coil (the transmitter). A second coil (the receiver) is positioned at a fixed distance from the first one, and the two coil axes may be aligned either vertically or horizontally with respect to the subsurface. The transmitting coil generates an electromagnetic (EM) field above the surface of the ground, a portion of which propagates into it.
This EM field, called the primary field $H_P$, induces an alternating electrical current within the ground, generating in turn a secondary EM field $H_S$, which propagates back to the surface and the air above. The second wire coil acts as a receiver, measuring the amplitude and phase components of the ratio between the primary and secondary EM fields. The complex measurements obtained by a GCM depend on some instrument settings, like the orientation of the dipoles, the frequency of the alternating current, the inter-coil distance, and the height of the instrument above the ground.

We provide below a brief discussion on the functional $\mathcal{M}$ presented in \eqref{eq:probl}. As we will see, the described model is one-dimensional, i.e., the measured data depends exclusively on the electrical conductivity below the measuring point. 

We assume that $N$ equispaced measurement sets $\mathbf{b}_j$ are performed in some interval $[a,b]$ (in our experiments we set $a=0$m and $b=10$m). Each set is obtained by collecting several measurements using different configurations of the FDEM device. Denoting a sampling of the electrical conductivity of the ground by
$$
\Sigmab=[\sigmab_1,\ldots,\sigmab_N],
$$
we can discretize the minimization problem \eqref{eq:probl} as follows
\begin{equation}\label{eq:model_decoupled}
\Sigmab=\arg\min_{\Sigmab}\sum_{j=1}^{N}\norm{M(\sigmab_j)-\mathbf{b}_j}_2^2,
\end{equation}
where $\norm{\cdot}_2$ is the Euclidean vector norm, $\boldsymbol{\Sigma}=[\sigmab_1,\ldots,\sigmab_N]\in\R^{n\times N}$, $\mathbf{B}=[\mathbf{b}_1,\ldots,\mathbf{b}_N]\in\C^{m\times N}$, and the vector function $\mathbf{M}(\Sigmab)=[M(\sigmab_1), \ldots,M(\sigmab_N)] \in \C^{m\times N}$ returns the readings predicted by the model in the same order they were arranged in the vector $\mathbf{b}_j$. This can be done because the model is one-dimensional, i.e., the values of the $j$th column of $\mathbf{B}$ depend only on the entries of the $j$th column of $\Sigmab$.

Therefore, once discretized, the two-dimensional problem is reduced to $N$ independent one-dimensional problems.

Although this may seem beneficial, since this decoupling can be greatly useful in lowering the computational complexity of the problem at hand, this is not completely true. A natural approach for solving the two-dimensional problem would be to solve each one-dimensional problem independently and then ``stack'' the obtained solutions in order to get a two-dimensional image; see \cite{ddrv19}. However, as we show in our numerical examples, in practice this approach may lead to poor reconstructions. In fact, the problem is severely ill-posed and in real-world applications the measured data is affected by noise, thus the obtained solution can deviate (also substantially) from the exact $\sigmab_j$ and this results in the stacked image to be ``spliced''. It is the purpose of this paper to propose a solution method that couples the $\sigmab_j$ in order to avoid this ``splicing'' and to counteract the ill-posedness of the problem. To this end, we consider the following variational problem 
\begin{equation}\label{eq:model}
	\arg\min_{\Sigmab\geq 0}\frac{1}{2}\norm{\mathbf{M}(\mathbf{\Sigma})-\mathbf{B}}_F^2+\frac{\gamma}{q}\norm{\mathbf{D}(\mathbf{\Sigma})}_q^q,
\end{equation}
where $\norm{\cdot}_F$ denotes the Frobenius norm, $\gamma>0$, and $0<q\leq 2$. With the notation $\mathbf{D}(\mathbf{\Sigma})$ we denote the vector in $\R^{nN}$ defined by
$$
\mathbf{D}(\mathbf{\Sigma})=D(\ts{vec}(\mathbf{\Sigma})),
$$
with
${D}={L}\otimes {I}+{I}\otimes {L}$, where $\otimes$ denotes the Kronecker product, ${I}$ is the identity matrix, $\ts{vec}$ is the vectorization operator, and 
$$
{L}=\begin{bmatrix}
2 &-1\\
-1&2&-1\\
  & \ddots &\ddots&\ddots\\
  &  &      -1&2&-1\\
  &  &      &  -1&2
\end{bmatrix}.
$$
The operator $\mathbf{D}$ (or equivalently the matrix $D$) is a discretization of the Laplacian, thus we are considering a high order version of so-called Total Variation regularization \cite{ROF}. When $q=2$ the model \eqref{eq:model} reduces to the well-known Tikhonov minimization for non-linear problems, here with the addition of the non-negativity constraint; see, e.g., \cite{TS56} for a first discussion on this kind of problems, \cite{Engl,Hansen} for a recent discussion on Tikhonov method and \cite{BBHRYZ17} for an algorithm for the solution of the non-negatively constrained Tikhonov minimization for linear problems. When $1<q<2$, the term $\norm{\mathbf{D}(\mathbf{\Sigma})}_q^q$ is convex and smooth, while for $q=1$ it is convex, but non-smooth. For $0<q<1$, we define $\norm{\mathbf{x}}_q^q=\sum_{j=1}^{n}|x_j|^q$ and we still refer to this quantity as to $\ell_q$-norm, however, this is not a norm since it does not satisfy the triangular inequality. Therefore, if $0<q<1$, the second term in \eqref{eq:model} is non-convex. Nevertheless, as shown in \cite{HLMRS17,BR19,BDDMR20,BPR20,BR20,BR21}, in imaging application it is beneficial to select $0<q<1$, this is due to the fact that, in this case, the $\ell_q$-norm approximates the $\ell_0$-norm, thus leading to reconstructions with sparse gradients. This property is desirable in an approximate reconstruction.

Finally, observe that, since we have used zero boundary conditions for the discretization of the Laplacian, it holds that $\ts{rank}({D})=nN$. Moreover,  note that, since we know that the electrical conductivity is positive, we imposed the non-negativity constraint.

The main contributions of this paper are the following. Firstly, we describe a new variational model for computing an approximate solution of \eqref{eq:probl}. The novelty of the proposed model is the horizontal coupling of the columns of $\Sigmab$ that allows us to obtain very accurate reconstructions of the electrical conductivity of the ground that are not affected by the ``splicing'' issue described above. To the best of our knowledge, this coupling has never been considered before in the literature for this kind of problem. Secondly, we prove the regularization properties of $\ell^2-\ell^q$ minimization for non-linear problems. This properties has been shown firstly in \cite{BPR20} for linear problems and we extend this analysis to the non-linear case. Finally, we propose an algorithm for the solution of the problem \eqref{eq:model} that extends to the non-linear case the one proposed in \cite{HLMRS17} that was constructed for linear problems. 

This paper is organized as follows: Section~\ref{sect:model} describes the functional $\mathcal{M}$ and its discretization, in Section~\ref{sect:reg} we show the theoretical properties of \eqref{eq:model}. Section~\ref{sect:min} provides an algorithm for the computation of an approximate solution of \eqref{eq:model}, in Section~\ref{sect:num} we show some numerical examples to demonstrate the performances of the proposed approach, and we draw some conclusions and outline future research in Section~\ref{sect:concl}.

\section{Modelization of the problem}\label{sect:model}
We now briefly describe how the FDEM device can be modeled. The main non-linear forward model which describes the interaction between the soil and the FDEM induction device when the electrical conductivity and the magentic permeability are known, has been described in~\cite{hendr02}.

In this model, the soil is assumed to have $n$ layers which are charaterized by an electrical conductivity $\sigma_i$ (measured in S/m) and a magnetic permeability $\mu_i$ (measured in H/m), for $i=1,\ldots,n$; see \cite{dfr14,ddr17} for more details. The thickness of each layer, measured in meters, is denoted by $d_i$, considering infinite the thikness of the deepest layer $d_n$. Finally, the distance between the coils is represented by $\rho$ and the height at which the measurements are taken by $h$.

Let  us now consider the propagation constant $u_i(\lambda) = \sqrt{\lambda^2 + \ii\sigma_i\mu_i\omega}$, where $\ii$ denotes the imaginary unit, ${\rm Re}(u_i(\lambda))\geq0$, and $\omega$ is the angular frequency of the instrument, that is, $2\pi$ times the frequency in Hertz. The variable $\lambda$ is the variable of integration which ranges from zero to infinity measuring the ratio between the depth in meters below the ground surface and the inter-coil distance $\rho$.

The surface admittance $Y_i(\lambda)$ at the top of each layer verifies the
recursion 
\begin{equation}
	\label{surfadm}
	Y_i(\lambda) = N_i(\lambda)\frac{Y_{i+1}(\lambda)+N_i(\lambda)
		\tanh(d_i u_i(\lambda))}{N_i(\lambda) + Y_{i+1}(\lambda)
		\tanh(d_i u_i(\lambda))},
\end{equation}
for $i=n-1,\ldots,1$, where $N_i(\lambda)= u_i(\lambda)/(\ii\mu_i\omega)$ represents the characteristic admittance at the $i$-th layer; see \cite{wait2}. Note that this formulation of $Y_i$ is not computationally stable and to solve practical problems different formulations have to be considered. However, our purpose here is to describe the model and we do not dwell further on the computational issues of the forward model. We refer the interested reader to \cite{ddflr20} for a discussion on the numerical implementation of the forward model.  At the last layer, i.e., when $i=n$, the characteristic admittance and the surface admittance coincide being $Y_n(\lambda) = N_n(\lambda)$ used to initialize the recursion \eqref{surfadm}.
We remark here, that both the characteristic and the surface admittances are functions
of the frequency $\omega$ via the functions $u_i(\lambda)$.

Now, the ratio of the secondary to the primary field for the vertical ($\nu=0$) and horizontal ($\nu=1$) orientation of the coils are given by
\begin{equation}\label{nonlinmodel}
		M_\nu(\sigmab,\mub;h,\omega,\rho) = -\rho^{3-\nu} \mathcal{H}_\nu\left[\lambda^{1-\nu} e^{-2h\lambda} R_{\omega,0}(\lambda) \right](\rho),\qquad \nu=0,1,
\end{equation}
where $\sigmab=(\sigma_1,\ldots,\sigma_n)^T$, $\mub=(\mu_1,\ldots,\mu_n)^T$,
$N_0(\lambda)=\lambda/(\ii\mu_0\omega)$, $\mu_0=4\pi\cdot 10^{-7}$H/m is the
magnetic permeability of free space. The reflection factor $R_{\omega,0}(\lambda)$ is defined by
\begin{equation*}
	R_{\omega,0}(\lambda) = \frac{N_0(\lambda) - Y_1(\lambda)}{N_0(\lambda) + Y_1(\lambda)},
\end{equation*}
with $Y_1(\lambda)$ computed by the recursion \eqref{surfadm}, and
$$
\mathcal{H}_\nu[f](\rho) = \int_{0}^{\infty} f(\lambda) J_\nu(\rho\lambda)
\lambda \,d\lambda , \qquad \nu=0,1,
$$
is the Hankel transform, where $J_0,J_1$ are first kind Bessel functions of order
0 and 1, respectively. 

Note that functions in \eqref{nonlinmodel} are complex valued functions. The imaginary part or quadrature component of the field ratio is usually interpreted as the apparent conductivity of the soil, while the magnetic permeability is related to the real or in-phase component.

Simultaneous measurements with different inter-coil distances or different operating frequencies can be recorded by recent FDEM devices at different heights.
We denote by $\rhob=(\rho_1,\ldots,\rho_{m_\rho})^T$, $\mathbf{h}=(h_1,\ldots,h_{m_h})^T$, and $\omegab=(\omega_1,\ldots,\omega_{m_\omega})^T$, the vectors containing the
loop-loop distances, the heights, and the angular frequencies at which the readings were taken. We consider the corresponding $m=2 m_\rho m_h m_\omega$ data points $b^\nu_{tls}$, where $t=1,\ldots,m_\rho$, $l=1,\ldots,m_h$, $s=1,\ldots,m_\omega$, while $\nu\in\{0, 1\}$ represents the vertical and horizontal orientations of the coils, respectively.
The observations $b^\nu_{tls}$ are rearranged in a vector $\mathbf{b}\in\C^m$.

In various papers this non-linear model has been studied for different device configurations and different techniques were applied; see, e.g., \cite{hendr02,dfr14,dr16,ddr17,ddrv18,ddrv19}. An algorithm for the regularized inversion of this model has been implemented in a Matlab package which includes a graphical user interface in \cite{ddflr20}.

We remark that, for small values of the conductivity of the soil, a linear model has been introduced in \cite{McNeill} and has been solved first in \cite{Borchers} and later in \cite{dfrv19} from the theoretical point of view, where an optimized solution method has been proposed.

In the following, it is assumed that the contribution of the permeability
distribution to the overall response is negligible, i.e. $\mu=\mu_0$, so that the measurements are considered to be sensitive merely to conductivity values. However, in principle, the regularization approach discussed here can be easily extended
to include also the inversion for the $\mub$ components; see \cite{ddr17}.
\section{Regularization property}\label{sect:reg}
In this section we consider the model \eqref{eq:model} in the presence of noise, that is

\begin{equation}\label{eq:modelnoise}
\arg\min_{\Sigmab\geq0}\frac{1}{2}\norm{\mathbf{M}(\mathbf{\Sigma})-\mathbf{B}^\delta}_F^2+\frac{\gamma}{q}\norm{\mathbf{D}(\mathbf{\Sigma})}_q^q,
\end{equation}
where $\delta$ denotes the noise present in the data, i.e.,
$$
\norm{\bm{B}-\bm{B}^\delta}\leq\delta.
$$
We would like to show that the model \eqref{eq:modelnoise} induces a regularization method. Namely, if $\gamma$ is chosen depending on the amount of noise that corrupts the data, then the solutions of \eqref{eq:modelnoise} converge to a minimum norm solution of the noise free problem \eqref{eq:model} as $\mathbf{B}^\delta\rightarrow\mathbf{B}$ or, equivalently, $\delta\rightarrow0$. Denote by $\mathcal{J}$ the functional minimized in \eqref{eq:modelnoise}, i.e.,
\begin{equation}\label{eq:J}
\mathcal{J}(\mathbf{\Sigma})=\frac{1}{2}\norm{\mathbf{M}(\mathbf{\Sigma})-\mathbf{B}^\delta}_F^2+\frac{\gamma}{q}\norm{\mathbf{D}(\mathbf{\Sigma})}_q^q+\iota_0(\Sigmab),
\end{equation}
where $\iota_0(\Sigmab)$ is the indicator function of the non-negative cone, namely
$$
\iota_0(\Sigmab)=\left\{\begin{array}{ll}
0&\mbox{if }\Sigmab\geq0,\\
\infty&\mbox{else.}
\end{array}\right.
$$

Before showing the regularization property we need to first prove some auxiliary results. The analysis of the model can be derived from the ones in \cite{BPR20,HKP07}, however, since some proofs are different, we report it here for the convenience of the reader.

We first recall the following result.
\begin{lem}[\cite{BPR20}]\label{lemma:bound}
Let $\{x_j\}_{j\in\N}$ be a sequence of elements of $\R^n$ and let $q>0$. If the $\norm{\bm{x}_j}_q^q$ are uniformly bounded, i.e., if there exists a constant $c>0$ independent of $j$ such that
\begin{equation*}
\norm{\bm{x}_j}_q^q\leq c\quad \forall j\in\N,
\end{equation*}
then $\norm{\bm{x}_j}_2^2$ is uniformly bounded.
\end{lem}

We can now show our first preliminary result.
\begin{prop}\label{prop:existence}
Let $\mathcal{J}$ be defined in \eqref{eq:J}, then $\mathcal{J}$ admits a global minimizer.
\end{prop}
\begin{prf}
	By definition $\mathcal{J}$ is lower semi-continuous and proper. Moreover, since $\ts{rank}(D)=nN$, it is easy to see that $\mathcal{J}$ is coercive. In fact, assume that the sequence $\{\Sigmab_j\}_{j\in\N}$ is such that $\norm{\Sigmab_j}_F^2\rightarrow\infty$ as $j\rightarrow\infty$, then, since $\ts{rank}(D)=nN$, $\norm{\mathbf{D}(\Sigmab_j)}_2^2\rightarrow\infty$ as $j\rightarrow\infty$. Consequently, $\norm{\mathbf{D}(\Sigmab_j)}_q^q\rightarrow\infty$ as $j\rightarrow\infty$ and thus $\mathcal{J}(\Sigmab_j)\rightarrow\infty$ as $j\rightarrow\infty$.
	
	Since $\mathcal{J}$ is proper, there exists $\Sigmab$ such that $\mathcal{J}(\Sigmab)<\infty$ and we can define
	$$
	\varphi=\inf_{\Sigmab\in\R^{n \times N}}\mathcal{J}(\Sigmab).
	$$
	By lower semi-continuity of $\mathcal{J}$, there exists a sequence $\{\Sigmab_j\}_{j\in\N}$ and $M\geq0$ such that 
	\begin{equation}\label{eq:conv_to_phi}
	\mathcal{J}(\Sigmab_j)\rightarrow\varphi \;\mbox{as}\; j\rightarrow\infty\quad\mbox{and}\quad\mathcal{J}(\Sigmab_j)\leq M\;\forall j\in\N.
	\end{equation}
	In particular, $\norm{\mathbf{D}(\Sigmab_j)}_q^q\leq M$ for all $j\in\N$ and, thanks to Lemma~\ref{lemma:bound}, there exists $\widetilde{M}$ such that $\norm{\mathbf{D}(\Sigmab_j)}_2^2\leq \widetilde{M}$ for all $j\in\N$. Since $\ts{rank}({D})=nN$, the sequence  $\{\Sigmab_j\}_{j\in\N}$ is uniformly bounded and, thus, it admits a convergent subsequence denoted by $\{\Sigmab_{j_k}\}_{j_k\in\N}$, with $\Sigmab_{j_k}\rightarrow\Sigmab^*$ as $j_k\rightarrow\infty$. We would like to show that $\Sigmab^*$ is a minimizer of $\mathcal{J}$. The definition of $\varphi$ yields
	\begin{equation*}
	\varphi\leq\mathcal{J}(\Sigmab^*)\leq\liminf_{j_k\rightarrow\infty}\mathcal{J}(\Sigmab_{j_k})=\lim_{j_k\rightarrow\infty}\mathcal{J}(\Sigmab_{j_k})=\varphi,
	\end{equation*}
	where the second to last equality follows from the lower semi-continuity of $\mathcal{J}$ and the last one follows from \eqref{eq:conv_to_phi}. This shows that $\mathcal{J}(\Sigmab^*)=\varphi$, i.e., that $\Sigmab^*$ is a global minimizer of $\mathcal{J}$ which concludes the proof.
\end{prf}
We are now in position to show our main result.
\begin{thm}
	Let $0<q\leq 2$ be fixed and $\{\mathbf{B}^{\delta_j}\}_{j\in\N}$ be a sequence such that, for all $j\in\N$, $\norm{\mathbf{B}-\mathbf{B}^{\delta_j}}_F\leq\delta_j$.
	Assume that $\delta_j\rightarrow0$ as $j\rightarrow\infty$. Let $\{\gamma_j\}_{j\in\N}$ be a sequence of positive real numbers such that
	$$
	\gamma_j\rightarrow0\quad\mbox{and}\quad\frac{\delta^2_j}{\gamma_j}\rightarrow0\quad\mbox{as}\quad j\rightarrow\infty.
	$$
	Denote by $\mathcal{J}_j$ the functional
	$$
	\mathcal{J}_j(\Sigmab)=\frac{1}{2}\norm{\mathbf{M}(\mathbf{\Sigma})-\mathbf{B}^{\delta_j}}_F^2+\frac{\gamma_j}{q}\norm{\mathbf{D}(\mathbf{\Sigma})}_q^q+\iota_0(\Sigmab)
	$$
	and, for all $j\in\N$, let
	$$
	\Sigmab_j\in\arg\min_{\Sigmab}\mathcal{J}_j(\Sigmab).
	$$
	Then there exists a convergent subsequence of $\{\Sigmab_j\}_{j\in\N}$, denoted by $\{\Sigmab_{j_k}\}_{j_k\in\N}$, such that
	$$
	\Sigmab_{j_k}\rightarrow\Sigmab^*\quad\mbox{as}\quad j_k\rightarrow\infty
	$$
	and
	$$
	\Sigmab^*\in\arg\min\{\norm{\mathbf{D}(\Sigmab)}_q^q:\;\mathbf{M}(\Sigmab)=\mathbf{B},\; \Sigmab\geq0\},
	$$
	assuming that this set is not empty.
\end{thm}
\begin{prf}
	First let us observe that, thanks to Proposition~\ref{prop:existence}, the sequence $\{\Sigmab_j\}_{j\in\N}$ is well defined.
	
	Since $\Sigmab_j$ is a global minimizer of $\mathcal{J}_j$, we have that for all $\Sigmab$ it holds
	$$
	\mathcal{J}_j(\Sigmab_j)\leq\mathcal{J}_j(\Sigmab).
	$$
	In particular, let $\Sigmab^\dagger\in\arg\min\{\norm{\mathbf{D}(\Sigmab)}_q^q:\;\mathbf{M}(\Sigmab)=\mathbf{B}\}$, then
	$$
	\mathcal{J}_j(\Sigmab_j)\leq\mathcal{J}_j(\Sigmab^\dagger).
	$$
	Recalling that $\mathbf{M}(\Sigmab^\dagger)=\mathbf{B}$ and that $\norm{\mathbf{B}-\mathbf{B}^{\delta_j}}_F\leq\delta_j$, for all $j\in\N$, we have that there exists $j_0$ such that, for all $j>j_0$, it holds
	\begin{align*}
		\frac{1}{2}\norm{\mathbf{M}(\mathbf{\Sigma}_j)-\mathbf{B}^{\delta_j}}_F^2+\frac{\gamma_j}{q}\norm{\mathbf{D}(\mathbf{\Sigma}_j)}_q^q&\leq\frac{1}{2}\norm{\mathbf{M}(\mathbf{\Sigma}^\dagger)-\mathbf{B}^{\delta_j}}_F^2+\frac{\gamma_j}{q}\norm{\mathbf{D}(\mathbf{\Sigma}^\dagger)}_q^q\nonumber\\
		&\leq\frac{1}{2}\delta_j^2+\frac{\gamma_j}{q}\norm{\mathbf{D}(\mathbf{\Sigma}^\dagger)}_q^q\nonumber\\
		&\leq C,\nonumber
	\end{align*}
	where $C$ is a constant independent of $j$, where the last inequality follows from the fact that $\delta_j,\gamma_j\rightarrow0$ as $j\rightarrow\infty$, and we observed that $\iota_0(\Sigmab_j)=\iota_0(\Sigmab^\dagger)=0$. In particular, we have that, for all $j>j_0$, $\norm{\mathbf{D}(\mathbf{\Sigma}_j)}_q^q\leq C$. Thanks to $\ts{rank}({D})=nN$ and Lemma~\ref{lemma:bound}, we have that $\norm{\mathbf{\Sigma}_j}_F\leq \widetilde{C}$ for a certain constant $\widetilde{C}$. Since the sequence $\{\Sigmab_j\}_{j\in\N}$ is uniformly bounded (for $j>j_0$), it admits a converging subsequence $\{\Sigmab_{j_k}\}_{j_k\in\N}$ and let $\Sigmab^*$ denote its limit. We first show that $\mathbf{M}(\Sigmab^*)=\mathbf{B}$.
	\begin{align*}
	0&\leq\frac{1}{2}\norm{\mathbf{M}(\Sigmab^*)-\mathbf{B}}^2_F\leq\liminf_{j_k\rightarrow\infty}\frac{1}{2}\norm{\mathbf{M}(\Sigmab_{j_k})-\mathbf{B}^{\delta_{j_k}}}^2_F\nonumber\\
	&\leq\liminf_{j_k\rightarrow\infty}\frac{1}{2}\norm{\mathbf{M}(\Sigmab_{j_k})-\mathbf{B}^{\delta_{j_k}}}^2_F+\frac{\gamma_{j_k}}{q}\norm{\mathbf{D}(\mathbf{\Sigma}_{j_k})}_q^q+\iota_0(\Sigmab_{j_k})\nonumber\\
	&\leq\liminf_{j_k\rightarrow\infty}\frac{\delta_{j_k}}{2}+\frac{\gamma_{j_k}}{q}\norm{\mathbf{D}(\mathbf{\Sigma}^\dagger)}_q^q+\iota_0(\Sigmab^\dagger)=0,\nonumber
	\end{align*} 
	i.e., $\mathbf{M}(\Sigmab^*)=\mathbf{B}$. We now show that $\Sigmab^*$ minimizes $\norm{\mathbf{D}(\mathbf{\Sigma}^*)}_q^q$. Recall that $\norm{\mathbf{D}(\mathbf{\Sigma}^\dagger)}_q^q$ is minimum, thus
	\begin{align*}
	\frac{1}{q}\norm{\mathbf{D}(\mathbf{\Sigma}^\dagger)}_q^q&\leq \frac{1}{q}\norm{\mathbf{D}(\mathbf{\Sigma}^*)}_q^q\leq\liminf_{j_k\rightarrow\infty}\frac{1}{q}\norm{\mathbf{D}(\mathbf{\Sigma}_{j_k})}_q^q+\iota_0(\Sigmab_{j_k})\nonumber\\
	&\leq\liminf_{j_k\rightarrow\infty}\frac{1}{2\gamma_{j_k}}\norm{\mathbf{M}(\Sigmab_{j_k})-\mathbf{B}^{\delta_{j_k}}}_F^2+\frac{1}{q}\norm{\mathbf{D}(\mathbf{\Sigma}_{j_k})}_q^q+\iota_0(\Sigmab_{j_k})\nonumber\\
	&=\liminf_{j_k\rightarrow\infty}\frac{1}{\gamma_{j_k}}\mathcal{J}_{j_k}(\Sigmab_{j_k})\\
	&\leq\liminf_{j_k\rightarrow\infty}\frac{1}{\gamma_{j_k}}\mathcal{J}_{j_k}(\Sigmab^\dagger)\nonumber\\
	&=\liminf_{j_k\rightarrow\infty}\frac{1}{2\gamma_{j_k}}\norm{\mathbf{M}(\Sigmab^\dagger)-\mathbf{B}^{\delta_{j_k}}}_F^2+\frac{1}{q}\norm{\mathbf{D}(\mathbf{\Sigma}^\dagger)}_q^q+\iota_0(\Sigmab^\dagger)\nonumber\\
	&\leq\liminf_{j_k\rightarrow\infty}\frac{\delta_{j_k}^2}{2\gamma_{j_k}}+\frac{1}{q}\norm{\mathbf{D}(\mathbf{\Sigma}^\dagger)}_q^q+\iota_0(\Sigmab^\dagger)=\frac{1}{q}\norm{\mathbf{D}(\mathbf{\Sigma}^\dagger)}_q^q.\nonumber
	\end{align*}
	We have shown that $\norm{\mathbf{D}(\mathbf{\Sigma}^*)}_q^q=\norm{\mathbf{D}(\mathbf{\Sigma}^\dagger)}_q^q$. Moreover, since the non-negative cone is a closed set it is obvious that $\Sigmab^*\geq0$, this coupled with the fact that $\mathbf{M}(\Sigmab^*)=\mathbf{B}$, concludes the proof.
\end{prf}

\section{Minimizing algorithm}\label{sect:min}
We wish now to discuss how to compute a solution of \eqref{eq:modelnoise}. Consider the minimization problem
\begin{equation}\label{eq:modelnoise_alt}
\arg\min_{\mathbf{\Sigma},\mathbf{\Xi}}\frac{1}{2}\norm{\mathbf{M}(\mathbf{\Sigma})-\mathbf{B}^\delta}_F^2+\frac{\gamma}{q}\norm{\mathbf{D}(\mathbf{\Xi})}_q^q+\frac{\beta}{2}\norm{\mathbf{\Sigma}-\mathbf{\Xi}}_F^2+\iota_0(\Sigmab),
\end{equation}
where $\beta>0$ is a fixed parameter. Obviously, if $\beta$ is large enough, the solutions of \eqref{eq:modelnoise} and \eqref{eq:modelnoise_alt} are the same. This reformulation is commonly performed in optimization so that the alternating minimization algorithm can be used. By applying alternating minimization to solve \eqref{eq:modelnoise_alt}, we obtain  the following iterations
\begin{equation}\label{eq:alt_min}
\left\{
\begin{array}{l}
\displaystyle \Sigmab^{(k+1)}=\arg\min_{\Sigmab}\frac{1}{2}\norm{\mathbf{M}(\mathbf{\Sigma})-\mathbf{B}^\delta}_F^2+\frac{\beta}{2}\norm{\mathbf{\Sigma}-\mathbf{\Xi}^{(k)}}_F^2+\iota_0(\Sigmab)\\ \displaystyle \mathbf{\Xi}^{(k+1)}=\arg\min_{\mathbf{\Xi}}\frac{\gamma}{q}\norm{\mathbf{D}(\mathbf{\Xi})}_q^q+\frac{\beta}{2}\norm{\mathbf{\Sigma}^{(k+1)}-\mathbf{\Xi}}_F^2.
\end{array}	\right.
\end{equation}
The convergence of the iterations in \eqref{eq:alt_min} is guaranteed by the results in \cite{GS99}. In particular, it holds the following
\begin{thm}
	Let $\mathbf{\Sigma}^{(k)}$ and $\mathbf{\Xi}^{(k)}$ denote the iterates defined in \eqref{eq:alt_min}, then
	$$
	(\mathbf{\Sigma}^{(k)},\mathbf{\Xi}^{(k)})\rightarrow(\mathbf{\Sigma}^{*},\mathbf{\Xi}^{*})\quad\mbox{as}\quad k\rightarrow\infty,
	$$
	where $(\mathbf{\Sigma}^{*},\mathbf{\Xi}^{*})$ is a stationary point of the functional in \eqref{eq:modelnoise_alt} that depends on the initial guess $\mathbf{\Xi}^{(0)}$. Moreover, if $\beta$ is large enough (possibly $\beta=\infty$), $\mathbf{\Sigma}^{*}=\mathbf{\Xi}^{*}=\widehat\Sigmab$, where $\widehat\Sigmab$ is a stationary point of the functional in \eqref{eq:modelnoise} that depends on $\mathbf{\Xi}^{(0)}$.
\end{thm}

\textbf{Remark.}Note that, since the minimized functional in \eqref{eq:modelnoise_alt} is non-convex it might have multiple global and local minima. Moreover, it may have several saddle points. However, in general, determining the global minimum of a non-convex function is a NP-hard problem. Therefore, in optimization theory, is usually considered a good enough result determining a stationary point of the minimized functional. Indeed, our numerical results show that the proposed approach is able to provide accurate reconstructions.

\subsection{Implementation details}
In the following, we detail how to numerically solve the two minimization subproblems in \eqref{eq:alt_min} at each iteration.

For the solution of the $\Sigmab$ subproblem in \eqref{eq:alt_min} let us first observe that we can write
$$
\mathbf{M}(\Sigmab)=[M(\sigmab_1),\ldots,M(\sigmab_N)],
$$
where $M(\sigmab_j)$ is a column vector. We can then rewrite the subproblem as
$$
\mathbf{\Sigma}^{(k+1)}=\arg\min_{\Sigmab}\sum_{j=1}^{N}\left[\frac{1}{2}\norm{M(\sigmab_j)-\mathbf{b}_j^\delta}_2^2+\frac{\beta}{2}\norm{\sigmab_j-\xib^{(k)}_j}_2^2+\iota_0(\sigmab_j)\right],
$$
where $\mathbf{b}_j^\delta$ and $\xib^{(k)}_j$ denote the $j$th column of $\mathbf{B}^\delta$ and $\mathbf{\Xi}^{(k)}$, respectively, and $\sigmab_j$ is the $j$th column of $\mathbf{\Sigma}$. Thus, we can write
$$
\sigmab_j^{(k+1)}=\arg\min_{\sigmab}\frac{1}{2}\norm{M(\sigmab)-\mathbf{b}_j^\delta}_2^2+\frac{\beta}{2}\norm{\sigmab-\xib^{(k)}_j}_2^2+\iota_0(\sigmab),\quad j=1,\ldots,N,
$$ 
i.e., the $\Sigmab$ subproblem decouples in $N$ independent one-dimensional subproblems. This problem can be rewritten as
\begin{equation}\label{eq:Sigma_sub}
\sigmab_j^{(k+1)}=\arg\min_{\sigmab}\frac{1}{2}\norm{\begin{bmatrix}M(\sigmab)\\\sqrt{\beta}\sigmab\end{bmatrix}-\begin{bmatrix}\mathbf{b}_j^\delta\\\sqrt{\beta}\xib^{(k)}_j\end{bmatrix}}_2^2+\iota_0(\sigmab),\quad j=1,\ldots,N.
\end{equation}
Observe that every $\sigmab_j^{(k+1)}$ subproblems are independent and thus can be solved in parallel. For the solution of the $\sigmab_j^{(k+1)}$, we consider a slight modification of the algorithm proposed in \cite{dfr14,dr16,ddr17,ddrv18,ddrv19,ddflr20}. Here we outline the algorithm and describe the modifications we made.

Let us first consider the following non-linear least-squares problem for a single column $M(\sigmab_j)$ and $\mathbf{b}_j^\delta$ of $\mathbf{M}(\Sigmab)$ and $\mathbf{B}^\delta$, respectively, i.e.
\begin{equation}\label{minimization}
\min_{\sigmab} \frac{1}{2} \norm{{M}(\sigmab)-\mathbf{b}_j^\delta}_2^2+\iota_0(\sigmab).
\end{equation}
We solve problem \eqref{minimization} by the Guass--Newton method.

We denote by $\mathbf{r}(\sigmab) = {M}(\sigmab)-\mathbf{b}^\delta_j$ the complex residual vector as a function of the conductivity $\sigmab$. At each step of the iterative algorithm we minimize the 2-norm of a linear approximation of the residual and compute $\sigmab^{(l+1)}=\sigmab^{(l)}+\mathbf{q}^{(l)}$, where
\begin{equation}\label{gaussnewt}
\displaystyle \mathbf{q}^{(l)}=\arg\min_{\mathbf{q}} \norm{\mathbf{r}(\sigmab^{(l)}) + J(\sigmab^{(l)})\mathbf{q}}_2,
\end{equation}
and $J(\sigmab^{(l)})$ is the Jacobian matrix of the function $\mathbf{r}$ computed in $\sigmab^{(l)}$.

Being the residual function $\mathbf{r}$ complex-valued, we solve 
problem \eqref{gaussnewt} by stacking the real and imaginary part of the
residual as follows (see \cite{ddrv19})
\begin{equation*}
\widetilde{\mathbf{r}}(\sigmab)= \begin{bmatrix}
\Real(\mathbf{r}(\sigmab)) \\
\Imag(\mathbf{r}(\sigmab))
\end{bmatrix} \in \R^{2m}, \quad
\widetilde{J}(\sigmab)=\begin{bmatrix}
\Real(J(\sigmab)) \\
\Imag(J(\sigmab))
\end{bmatrix} \in \R^{2m\times n}.
\end{equation*}
In the same way we set $\widetilde{\mathbf{r}}$ and $\widetilde J$, we rearrange the vectors $\mathbf{M}$ and $\mathbf{b}^\delta$, and we denote them by  $\mathbf{\widetilde M}$ and $\widetilde{\mathbf{b}^\delta}$, respectively. So, we replace \eqref{gaussnewt} by
\begin{equation}
\min_{\mathbf{q}}
\norm{\widetilde{\mathbf{r}}(\sigmab^{(l)})+\widetilde{J}_l \mathbf{q}}_2,
\label{leastsquaresnew}
\end{equation}
with $\widetilde{J}_l=\widetilde{J}(\sigmab^{(l)})$,
and the iterative method becomes
$$
\sigmab^{(l+1)} = \sigmab^{(l)} + \alpha_l\mathbf{q}^{(l)} 
= \sigmab^{(l)} - \alpha_l \widetilde{J}_l^{\dagger} \, \widetilde{\mathbf{r}}(\sigmab^{(l)}),
$$
where $\widetilde{J}_l^{\dagger}$ is the Moore--Penrose pseudoinverse of $\widetilde{J}_l$ and
$\alpha_l$ is a damping parameter which ensures the convergence.
This parameter is determined by coupling the Armijo--Goldstein principle~\cite{bjo96} to
the positivity constraint $\sigmab^{(l+1)}\geq0$ (see \cite{ddr17,dfr14}) so that
\begin{equation}\label{eq:alphal}
\norm{\mathbf{\widetilde r}(\mathbf{\sigma}^{(l)})}_2^2 - \norm{\mathbf{\widetilde r}(\mathbf{\sigma}^{(l)} + \alpha_l \mathbf{q}^{(l)})}_2^2 \geq \frac{1}{2} \alpha_l \norm{\widetilde J_l\mathbf{q}^{(l)}}_2^2\quad\mbox{with}\quad\sigmab^{(l+1)}\geq0
\end{equation}
is verified. The analytical expression of the Jacobian matrices with respect to the electrical conductivity and the magnetic permeability were computed in~\cite{dfr14} and~\cite{ddr17}, respectively, where it has also been proved that the computation of the analitical expression of the Jacobian matrices is faster than their finite difference approximations.

It is well known that the minimization problem~\eqref{leastsquaresnew} is
extremely ill-posed, meaning that the matrix $\widetilde{J}_l$ is severely ill-conditioned. In order to overcome this difficulty, we apply the generalized truncated singular value decomposition (GTSVD), as in \cite{dfr14,dr16,ddr17,ddrv19,ddflr20}, for stably computing an approximation of $\widetilde{J}_l^{\dagger}$. We introduce
a regularization matrix $\widehat{L}\in\R^{p\times 2n}$ ($p\leq 2n$), whose null space approximately contains the sought solution~\cite{rr13}. Under the assumption $\mathcal{N}(\widetilde{J}_l)\cap\mathcal{N}(\widehat{L})=\{\mathbf{0}\}$,
problem~\eqref{leastsquaresnew} is replaced by
\begin{equation}\label{eq:q}
\min_{\mathbf{q}\in\mathcal{S}} \norm{L\mathbf{q}}_2^2, \qquad
\mathcal{S} = \{ \mathbf{q}\in\R^{2n} \,:\, \widetilde{J}_l^T\widetilde{J}_l\mathbf{q}=-\widetilde{J}_l^T\widetilde{\mathbf{r}}(\sigmab^{(l)}) \}.
\end{equation}
Note that, thanks to the assumption above on the null spaces of $\widetilde{J}_l$ and $\widehat{L}$, the solution of \eqref{eq:q} is unique. Very common choices for $L$ are the discretization of the first or
second derivative operators.

Let the generalized singular value decomposition (GSVD)~\cite{gvl96} of the
matrix pair $(\widetilde{J}_l,\widehat{L})$ be
$$
\widetilde{J}_l = U \Sigma_J Z^{-1}, \qquad \widehat{L} = V \Sigma_L Z^{-1},
$$
where $U$ and $V$ are matrices with orthonormal columns $\mathbf{u}_i$ and $\mathbf{v}_i$, respectively, $Z$ is a non-singular matrix with columns $\mathbf{z}_i$, and $\Sigma_J$, $\Sigma_L$ are diagonal matrices with diagonal entries $c_i$ and $s_i$, which are the singular values of $\widetilde{J}$ and $\widehat{L}$, respectively.
Under the assumption that $m=2 m_\rho m_h m_\omega<2n$, common in the generality of cases, the truncated GSVD (TGSVD) solution $\mathbf{q}^{(l)}$ (see~\cite{Hansen} for details) can be written as 
\begin{equation}\label{eq:qGSVD}
\mathbf{q}^{(l)}_{\rm GSVD} = -\sum_{i=p-\ell+1}^p
    \frac{\mathbf{u}_i^T\widetilde{\mathbf{r}}^{(l)}}{c_{i-2n+\kappa}}\, \mathbf{z}_i
    - \sum_{i=p+1}^{2n} (\mathbf{u}_i^T\widetilde{\mathbf{r}}^{(l)})\, \mathbf{z}_i,
\end{equation}
where $\ell$ is fixed, $\kappa=\rank(\widetilde{J}_l)$, $\ell=1,\ldots,\kappa+p-2n$
is the regularization parameter, and $\mathbf{r}^{(l)}=\mathbf{r}(\sigmab^{(l)})$.

The resulting regularized damped Gauss--Newton method reads
$$
\sigmab^{(l+1)} =
\sigmab^{(l)} + \alpha_l \mathbf{q}^{(l)}_{\rm GSVD},
$$
with $\ell$ fixed and $\alpha_l$ determined at each step as in \eqref{eq:alphal}.

We now apply this algorithm to our specific case. In particular, we have to solve \eqref{eq:Sigma_sub}. We apply the previously described procedure to the following function getting
\begin{equation}\label{jacobian}
\widehat{M}(\sigmab)= \begin{bmatrix}
\widetilde M(\sigmab)\\\sqrt{\beta}\sigmab
\end{bmatrix}, \quad
\widehat{J}(\sigmab)=\begin{bmatrix}
\widetilde J(\sigmab)\\\sqrt{\beta}I
\end{bmatrix},
\end{equation}
where $I$ is the identity matrix. We remind here that $\widehat{M}$ is not a complex valued function anymore since we already stacked the real and imaginary parts of $M$ and $J$.

Let $\sigmab_j^{(k,l)}$ be an approximation of $\sigmab_j^{(k+1)}$. Denote by $\widehat{J}_l$ the Jacobian of $\widehat{M}$ in $\sigmab_j^{(k,l)}$ and consider the GSVD of the pair $\{\widehat{J}_l,\widehat L\}$ 
$$
\widehat{J}_l = \widehat U \widehat \Sigma_J \widehat Z^{-1}, \qquad \widehat L = \widehat V \widehat \Sigma_L \widehat Z^{-1}.
$$ 
Let 
$$
\widehat{\bm{r}}^{(l)}=\widehat{M}(\sigmab_j^{(k,l)})-\begin{bmatrix}
\bm{b}_j^\delta\\\xib^{(k)}_j
\end{bmatrix}=\begin{bmatrix}
M(\sigmab_j^{(k,l)})-\bm{b}_j^\delta\\\sqrt{\beta}(\sigmab_j^{(k,l)}-\xib^{(k)}_j)\end{bmatrix},
$$
then we compute, analogously to \eqref{eq:qGSVD},
$$
\mathbf{q}^{(l)}_{\rm GSVD} = -\sum_{i=p-\ell+1}^p
\frac{\mathbf{u}_i^T\widehat{\mathbf{r}}^{(l)}}{c_{i-2n+\kappa}}\, \mathbf{z}_i
- \sum_{i=p+1}^{2n} (\mathbf{u}_i^T\widehat{\mathbf{r}}^{(l)})\, \mathbf{z}_i,
$$
leading to the iteration
$$
\sigmab_j^{(k+1,l)}=\sigmab_j^{(k,l)}+\alpha_l\mathbf{q}^{(l)}_{\rm GSVD},
$$
where $\alpha_l$ is determined by the Armijo-Goldstein rule so that
$$
\norm{\widehat{\mathbf{r}}^{(l)}}^2 - \norm{\widehat{\mathbf{r}}^{(l+1)}}_2^2 \geq \frac{1}{2} \alpha_l \norm{\widehat{J}_l\mathbf{q}^{(l)}_{\rm GSVD}}_2^2\quad\mbox{and}\quad\sigmab^{(k,l+1)}\geq0.
$$

We now move to the $\mathbf{\Xi}$ subproblem. To solve this problem we consider the majorization-minimization algorithm proposed by Huang et al. in \cite{HLMRS17} and furtherly developed in \cite{BDDMR20,BR19,BR20}. We briefly describe the algorithm presented in \cite{HLMRS17}. To simplify the computations we consider a modified version of the $\mathbf{D}$ operator, namely we impose reflexive boundary conditions to the discretization of the Laplacian. In this way the obtained matrix has an exploitable structure that helps in the computations. In detail, with abuse of notation, we write 
$$
{L}=\begin{bmatrix}
1&-1\\
-1&2&-1\\
&  \ddots&\ddots&\ddots\\
&  &      -1&2&-1\\
&  &      &  -1&1
\end{bmatrix}.
$$
Thanks to the structure of $L$ the matrix $D$, which is defined by $D=L\otimes I+I\otimes L$, is the sum of a block Hankel with Hankel blocks matrix, a block Toeplitz with Hankel blocks matrix, a block Hankel with Toeplitz blocks matrix, and a block Toeplitz with Toeplitz blocks matrix. We recall that Toeplitz matrices are matrices that are constant on the diagonals and Hankel matrices are matrices that are constant on the anti-diagonals. Since $[-1\;2\;-1]$ is symmetric, if ${C}$ denotes the discrete cosine transform matrix, we have that
\begin{equation}\label{eq:fact}
{D}={C}^T{\Lambda}{C},
\end{equation}
where ${\Lambda}$ is a diagonal matrix. The diagonal elements of ${\Lambda}$ are computed as the cosine coefficients of the first column of ${D}$; see \cite{NC99} for more details.

We now describe the Majorization-Minimization (MM) procedure that we use for solving the $\mathbf{\Xi}$ subproblem. Let us first rewrite the minimization problem in \eqref{eq:alt_min} as follows
\begin{equation}\label{eq:2probl_simpl}
\ts{vec}(\mathbf{\Xi}^{(k+1)})=\arg\min_{\xib}\frac{1}{2}\norm{\xib-\sigmab^{(k+1)}}_2^2+\frac{\gamma}{q\beta}\norm{{D}\xib}_q^q,
\end{equation}
where $\sigmab^{(k+1)}=\ts{vec}(\Sigmab^{(k+1)})$. The MM algorithm generates a sequence of vectors $\xib^{(k,l)}$ that converges to an approximate solution of \eqref{eq:2probl_simpl}. Firstly, if $q\leq 1$, we need to smooth the $q$-norm so that it is differentiable. Let $\varepsilon>0$ be a small constant, then, for $\mathbf{x}\in\R^n$, it holds
$$
\norm{\mathbf{x}}_q^q=\sum_{j=1}^n |x_j|^q\approx\sum_{j=1}^n \left(x_j^2+\varepsilon^2\right)^{q/2}=:\norm{\mathbf{x}}_{q,\varepsilon}^q.
$$
Note that the function $\mathbf{x}\mapsto\norm{\mathbf{x}}_{q,\varepsilon}^q$ is differentiable everywhere. We consider the smoothed problem
\begin{equation}\label{eq:2probl_smooth}
\ts{vec}(\mathbf\Xi^{(k+1)})=\arg\min_{\xib}\frac{1}{2}\norm{\xib-{\sigmab}^{(k+1)}}_2^2+\frac{\gamma}{q\beta}\norm{{D}\xib}_{q,\varepsilon}^q=\arg\min_{\xib}\widehat{\mathcal{J}}_\varepsilon(\mathbf{\xib}).
\end{equation}
As pointed out in \cite{BR19}, the solutions of \eqref{eq:2probl_simpl} and \eqref{eq:2probl_smooth} are extremely similar and this substitution does not have any negative effect.

Let $\xib^{(k,l)}$ be an approximate solution of \eqref{eq:2probl_smooth}. We first construct a quadratic tangent majorant of $\widehat{\mathcal{J}}_\varepsilon$ in $\xib^{(k,l)}$ that majorizes it, i.e., a quadratic function $\mathcal{Q}(\xib,\xib^{(k,l)})$ such that
\begin{itemize}
	\item $\mathcal{Q}(\xib,\xib^{(k,l)})\geq \widehat{\mathcal{J}}_\varepsilon(\xib)$ for all $\xib$;
	\item $\mathcal{Q}(\xib^{(k,l)},\xib^{(k,l)})= \widehat{\mathcal{J}}_\varepsilon(\xib^{(k,l)})$;
	\item$\nabla\mathcal{Q}(\xib^{(k,l)},\xib^{(k,l)})= \nabla\widehat{\mathcal{J}}_\varepsilon(\xib^{(k,l)})$.
\end{itemize}
In \cite{HLMRS17} the authors provide two different choices for the construction of such a functional, here we considered the so-called fixed approach. Denote by $\widetilde{\mathbf{u}}^{(l)}=D\xib^{(k,l)}$, then we compute the vector $\mathbf{u}^{(l)}$ as
$$
\mathbf{u}^{(l)}=\widetilde{\mathbf{u}}^{(l)}\left(1-\left(\frac{(\widetilde{\mathbf{u}}^{(l)})^2+\varepsilon^2}{\varepsilon^2}\right)^{q/2-1}\right),
$$
where all the operations are meant element-wise. It is possible to see that the function 
$$
\mathcal{Q}(\xib,\xib^{(k,l)})=\frac{1}{2}\norm{\xib-\sigmab^{(k+1)}}_2^2+\frac{\gamma\varepsilon^{q-2}}{2\beta}\left(\norm{{D}\xib}_2^2-2\langle\mathbf{u}^{(l)},{D}\xib\rangle\right)+c,
$$
where $c$ is a constant independent of $\xib$, is a quadratic tangent majorant; see \cite{HLMRS17} for a derivation. An improved approximation of the solution of \eqref{eq:2probl_smooth} can be obtained as the unique minimizer of $\mathcal{Q}(\xib,\xib^{(k,l)})$, i.e., 
\begin{align*}
\xib^{(k,l+1)}&=\arg\min_{\xib}\frac{1}{2}\norm{\xib-{\sigmab}^{(k+1)}}_2^2+\frac{\gamma\varepsilon^{q-2}}{2\beta}\left(\norm{D\xib}_2^2-2\langle\mathbf{u}^{(l)},{D}\xib\rangle\right)\\
&=\arg\min_{\xib}\frac{1}{2}\norm{\begin{bmatrix}
\mathbf{I}\\\sqrt{\eta}{D}
\end{bmatrix}\xib-\begin{bmatrix} {\sigmab^{(k+1)}}\\\sqrt{\eta}\mathbf{u}^{(l)}\end{bmatrix}}_2^2,
\end{align*}
where $\eta=\frac{\gamma\varepsilon^{q-2}}{\beta}$. Writing the normal equation yields
$$
({I}+\eta{D}^T{D})\xib^{(k,l+1)}=\sigmab^{(k+1)}+\eta D^T\mathbf{u}^{(l)}.
$$
Using the factorization \eqref{eq:fact} we get
$$
{C}^T\left({I}+\eta{\Lambda}^T{\Lambda}T\right){C}\xib^{(k,l+1)}=\sigmab^{(k+1)}+\eta {D}^T\mathbf{u}^{(l)}.
$$
Thus, we can compute $\xib^{(k,l+1)}$ as
$$
\xib^{(k,l+1)}={C}^T\left({I}+\eta\Lambda^T\Lambda\right)^{-1}{C}(\sigmab^{(k+1)}+\eta {D}^T\mathbf{u}^{(l)}),
$$
where the inversion is well-defined since $\eta>0$. These computations can be performed fairly inexpensively, since the system to solve is a diagonal one and the application of the $n\times n$ cosine matrix can be performed in $O(n\log n)$ operations by means of the \texttt{dct} algorithm.

We summarize all the computations in Algorithm~\ref{algo:final}.

\begin{algorithm}
	\SetKwInOut{Input}{input}
	\SetKwProg{Parfor}{parfor}{}{end}
		\Input{$\mathbf{B}^\delta$, $\mathbf{\Xi}^{(0)}$, $\gamma>0$, $\beta>0$, $0<q\leq2$, $\ell\ll nN$, $\varepsilon>0$, $\sigma_0>0$.}
		$\eta=\frac{\gamma\varepsilon^{q-2}}{\beta}$\;
		\For{$k=0,1,2,\ldots$}{
		\Parfor{$j=1,\ldots,N$}{
		$\sigmab^{(k,1)}=\sigma_0\bm{1}$\;
		\For{$l=1,2,\ldots$}{
		$\widehat{J}_l=\begin{bmatrix}
				{\widetilde{J}}(\sigmab^{(k,l)}) \\
				\sqrt{\beta} I
				\end{bmatrix}$\;
		Compute the GSVD ${\widehat J} = \widehat U \widehat{\Sigma}_J \widehat Z^{-1}, \; \widehat{L} = \widehat V \widehat{\Sigma}_L \widehat Z^{-1}$\;
		${\bm{r}}^{(l)}=\begin{bmatrix}
			M(\sigmab_j^{(k,l)})-\bm{b}_j^\delta\\\sqrt{\beta}(\sigmab_j^{(k,l)}-\xib^{(k)}_j)\end{bmatrix}$\;
		$\displaystyle\mathbf{q}^{(l)}_{\rm GSVD} = -\sum_{i=p-\ell+1}^p \frac{\mathbf{u}_i^T\mathbf{r}^{(l)} }{c_{i-2n+\kappa}}\, \mathbf{z}_i - \sum_{i=p+1}^{2n} (\mathbf{u}_i^T\mathbf{r}^{(l)})\, \mathbf{z}_i$\;
		Find $\alpha_l$ s.t. 
		$
		\left\{ \begin{array}{l} \norm{\mathbf{r}^{(l)}}^2 - \norm{\begin{bmatrix}
		M\left(\mathbf{\sigma}^{(k,l)} + \alpha_l \mathbf{q}_{\rm GSVD}^{(l)}\right)-\bm{b}_j^\delta\\\sqrt{\beta}(\boldsymbol{\sigma}^{(k,l)} + \alpha_l \mathbf{q}_{\rm GSVD}^{(l)}-\xib^{(k)}_j)\end{bmatrix}}^2 \geq \frac{1}{2} \alpha_l \norm{\widehat{J}_l\mathbf{q}_{\rm GSVD}^{(l)}}^2 \\ \boldsymbol{\sigma}^{(k,l)} + \alpha_l \mathbf{q}_{\rm GSVD}^{(l)}\geq0\end{array}\right.
		$\;
		$\boldsymbol{\sigma}^{(k,l+1)}=\boldsymbol{\sigma}^{(k,l)} + \alpha_l \mathbf{q}_{\rm GSVD}^{(l)}$\;
	}
		$\boldsymbol{\sigma}^{(k+1)}=\boldsymbol{\sigma}^{(k,\infty)}$
	}

		$\sigmab^{(k+1)}=\ts{vec}(\Sigmab^{(k+1)})$\;
		$\xib^{(k,1)}=\sigmab^{(k+1)}$\;
		\For{$l=1,2,\ldots$}{
		$\widetilde{\mathbf{u}}^{(l)}=D\xib^{(k,l)}$\;
		 $\mathbf{u}^{(l)}=\widetilde{\mathbf{u}}^{(l)}\left(1-\left(\frac{(\widetilde{\mathbf{u}}^{(j)})^2+\varepsilon^2}{\varepsilon^2}\right)^{q/2-1}\right)$ \;
 	$\xib^{(k,l+1)}={C}^T\left({I}+\eta\Lambda^T\Lambda\right)^{-1}{C}(\sigmab^{(k+1)}+\eta{D}^T\mathbf{u}^{(j)})$\;}	
		$\ts{vec}(\mathbf\Xi^{(k+1)})=\xib^{(k,\infty)}$\;
	}

	\caption{Parallel alternating minimization for the solution of \eqref{eq:modelnoise_alt}}
	\label{algo:final}
\end{algorithm}
\section{Numerical examples}\label{sect:num}
In the numerical tests illustrated in this section we consider both synthetic and experimental data. The two-dimensional representation of the electrical conductivity is determined under the assumption that the magnetic permeability is the one of the free space, i.e., $\mub\equiv\mu_0$. 

Throughout this section, we show the effectiveness of our method and we compare the results obtained from the algorithm presented in this paper with those obtained by applying the method described in \cite{ddrv19}, i.e., by solving each one-dimensional problem independently and then ``stacking'' the obtained solutions.

All the computations are performed on an Intel(R) Xeon(R) Gold 6136 CPU @ 3.00GHz computer with 128Gb of RAM memory and 32 cores, running the Debian GNU/Linux operating system and Matlab R2020b.
\subsection{Synthetic data}
In the first test we consider two different configurations of the FDEM device. We show that the proposed algorithm is able to reconstruct an accurate approximation of the electrical conductivity for both devices. In out second test we report the solution obtained by varying the dimensions of the problem discretization for one configuration.

We generate the data matrix $\bm{B}_{\rm{exact}}$ of a chosen dimension $m \times N$ and we simulate the presence of noise in the data, by letting
$$
\bm{e} =\frac{\delta}{\sqrt{m}} \norm{\bm{B}_{\rm{exact}}}_F^2 \bm{w},
$$
where $\bm{w}$ is a vector with normally distributed entries having zero mean and unit variance, and $\delta$ represents the noise level. The new data matrix is denoted by $\bm{B}^\delta$.
\paragraph{Test 1.}
This example concerns the reconstruction of a two-dimensional model generated from $50$ soundings along a $10$ m straight line and $20$ layers in detph. It is characterized by an increasing change in conductivity (from $0$ S/m to $1$ S/m) occurring at an increasing depth. 

The synthetic data simulate an acquisition performed by both the Geophex GEM-2 and CMD Explorer instruments, with two orientations of the coils and one measurement height $h = 1$ m, i.e., the instrument is kept $1$ m above the ground. We consider a slice of the ground that is $10$ meters wide. The first device, the Geophex GEM-2, works with an intercoil distance of $1.66~{\rm m}$ ($\rho = 1.66~{\rm m}$) and six different frequencies $f = 775, 1175, 3925, 9825, 21725, 47025~{\rm Hz}$, while the second one, the CDM Explorer, measures with three different intercoil distances $\rho = 1.48, 2.82, 4.49~{\rm m}$ and only one frequency $f = 10~{\rm kHz}$. 
The data values are finally perturbed by uncorrelated Gaussian noise with standard deviation $\delta = 10^{-2}$. Here, the choice of the regularization parameter $\ell$ does not affect substantially the computed results. This is due to the structure of the Jacobian matrix defined in \eqref{jacobian}. The matrix $\widehat{J}$ is obtained by stacking the Jacobian of ${M}$ and a scaled identity matrix, therefore, the number of the singular values does not change, but they do not decrease as rapidly as the one of $\widetilde{J}$. Therefore, we set in all our experiments $\ell=15$. From our experience we observed that a larger value of $\ell$ does not neither improve nor deteriorate the obtained approximate solutions. In the algorithm described in \cite{ddrv19} the Jacobian matrix of the problem was very ill-conditioned and the results depended on the choice of the truncation parameter. We have established a maximum number of $50$ iterations, $p=2$, $q=0.1$, the regularization parameter $\gamma = 10^{-4}$, and $\sigma_0=0.1$.

\begin{figure}
\centering
\begin{minipage}{0.4\linewidth}
	\centering
	\includegraphics[width=\linewidth]{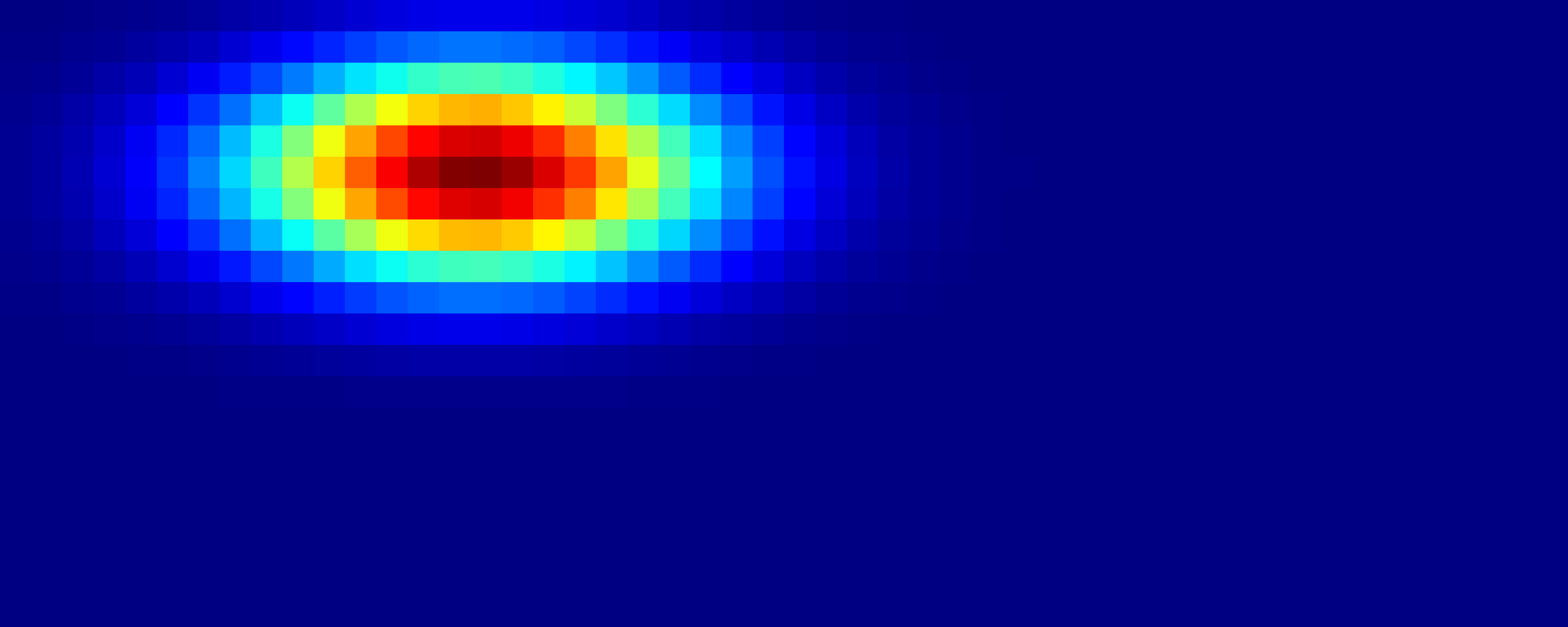}\\(a)
\end{minipage}
\begin{minipage}{0.4\linewidth}
	\centering
	\includegraphics[width=\linewidth]{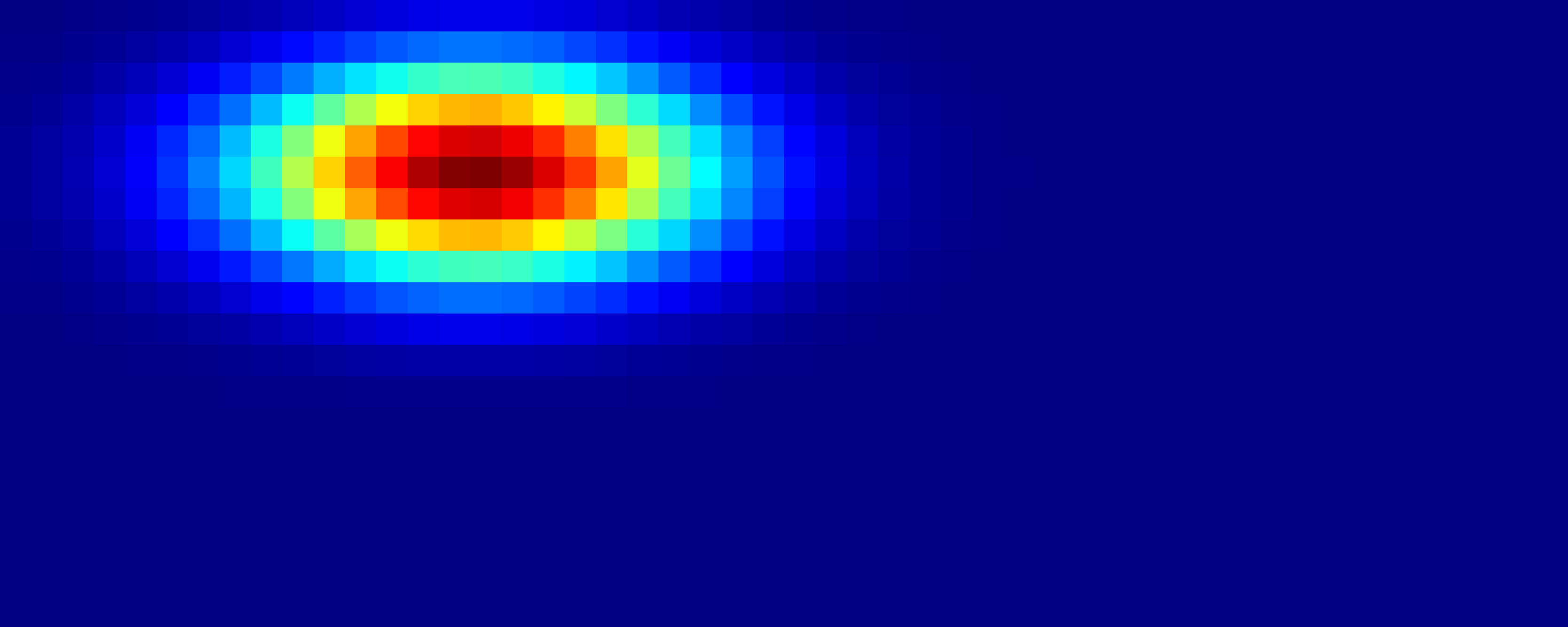}\\(b)
\end{minipage}

\begin{minipage}{0.4\linewidth}
	\centering
	\includegraphics[width=\linewidth]{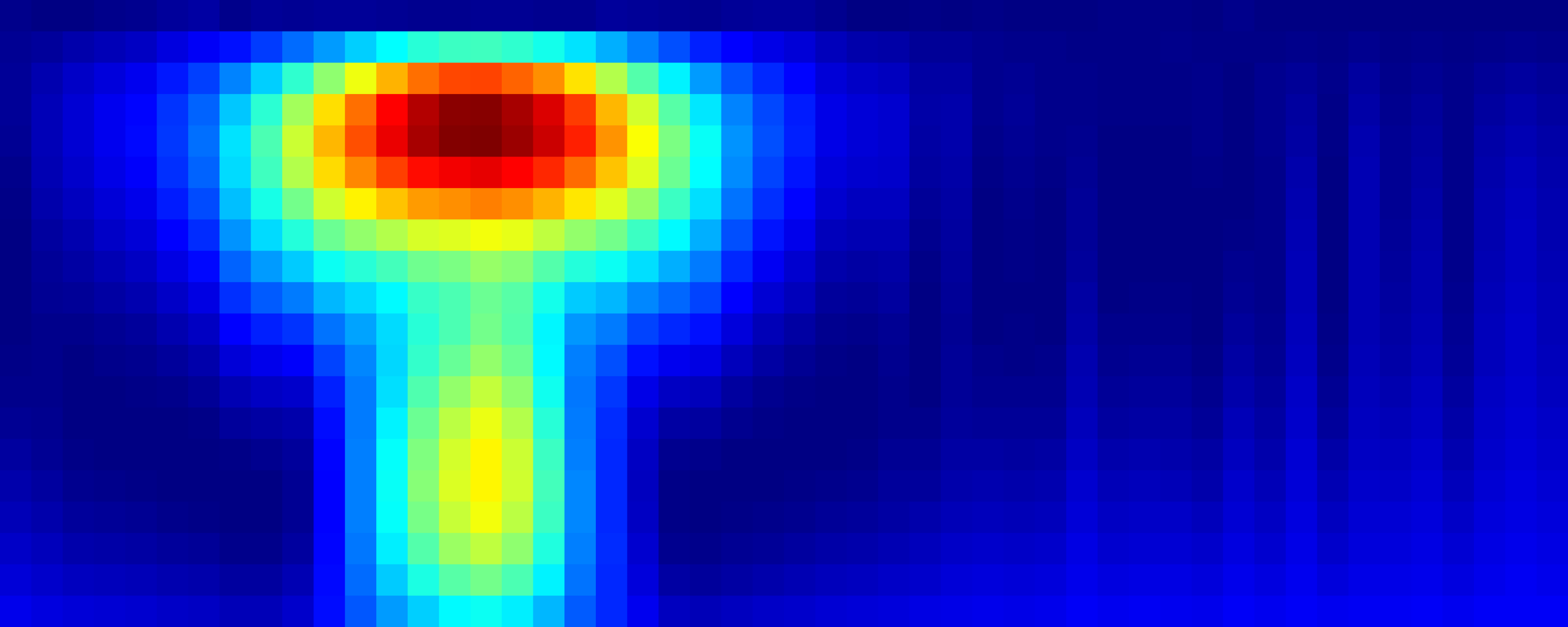}\\(c)
\end{minipage}
\begin{minipage}{0.4\linewidth}
	\centering
	\includegraphics[width=\linewidth]{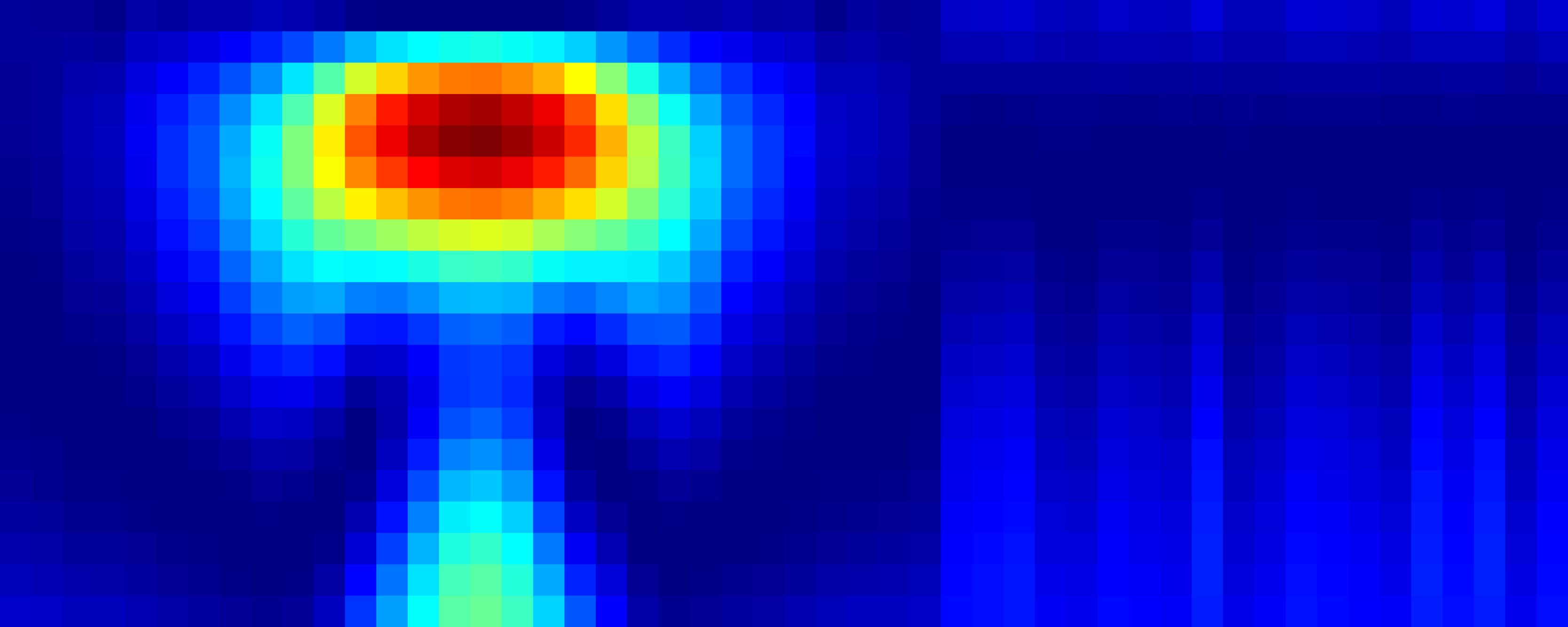}\\(d)
\end{minipage}

\begin{minipage}{0.4\linewidth}
	\centering
	\includegraphics[width=\linewidth]{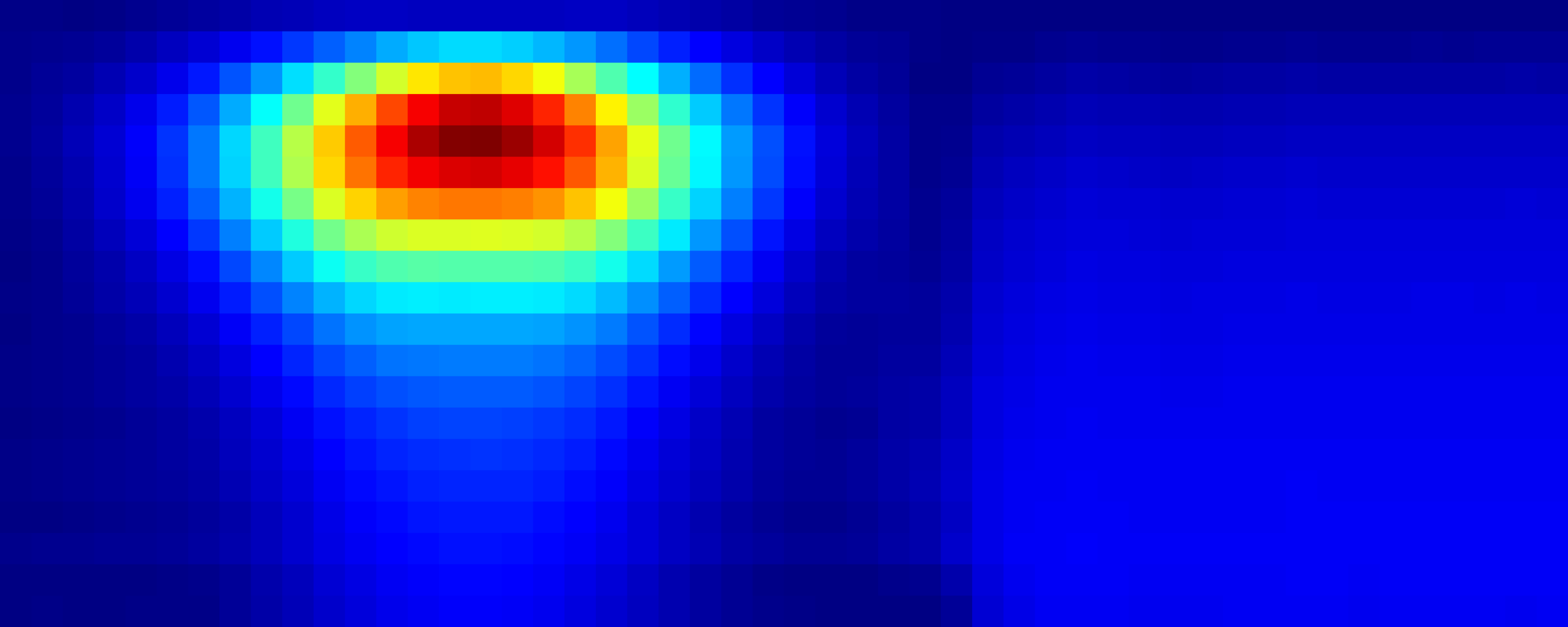}\\(e)
\end{minipage}
\begin{minipage}{0.4\linewidth}
	\centering
	\includegraphics[width=\linewidth]{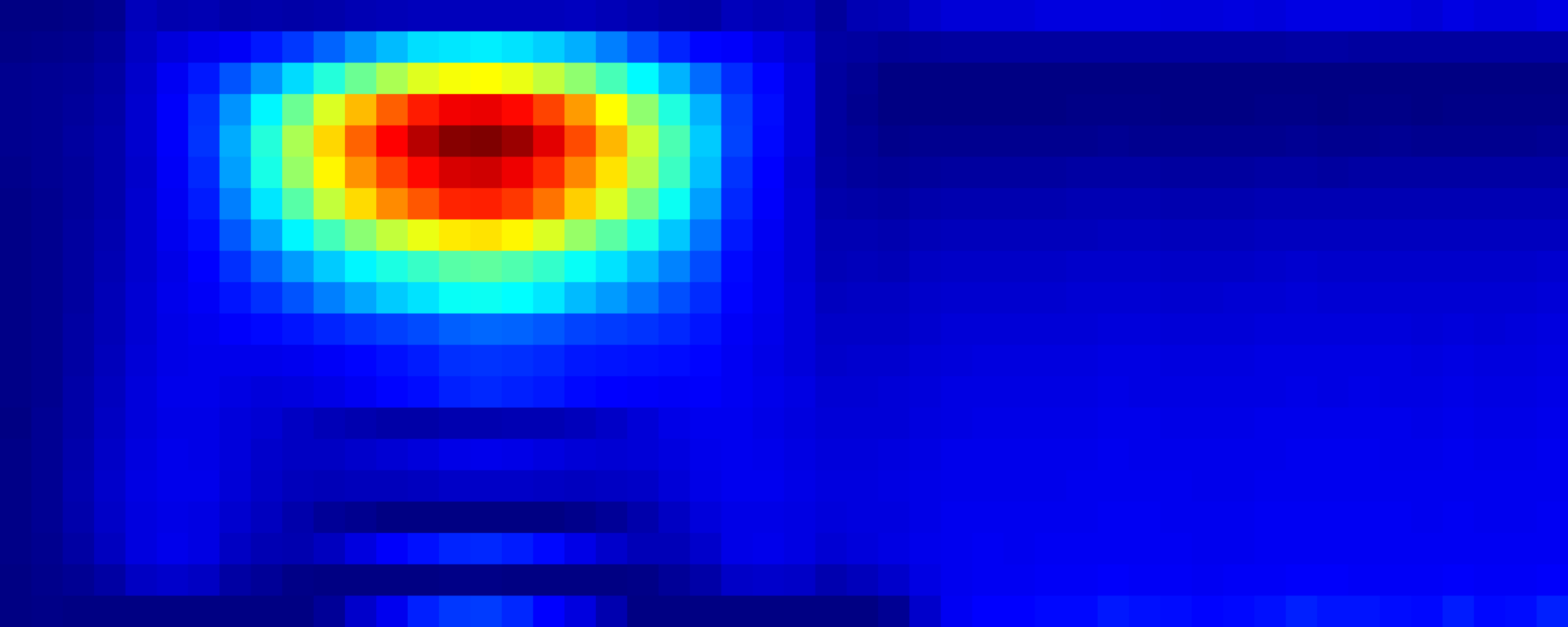}\\(f)
\end{minipage}
\caption{Reconstruction of the electrical conductivity from two devices configuration. The first row represents the exact solution used to generate the data with both the Geophex - GEM 2 (panel (a)) and the CDM Explorer (panel (b)). On the second row, is reported the reconstruction of the one-dimensional models side by side in a pseudo two-dimensional fashion (panels (c) and (d)). On the last row, the obtained results by the algorithm presented in this paper are depicted (panels (e) and (f)). The data values are corrupted by Gaussian noise with noise level $\delta = 10^{-2}$, the maximum number of iterations is fixed at $50$, $p=2$, $q=0.1$, and the regularization parameter $\gamma = 10^{-4}$.}
\label{Figure2}
\end{figure}

In the first row of Figure~\ref{Figure2}, we report the exact solution used to generate the data using both the devices: Geophex - GEM 2 on the left, CDM Explorer on the right. On the middle row, we report the reconstruction of the one-dimensional models side by side in a pseudo two-dimensional fashion as in \cite{ddrv19}. As we described above, the method in \cite{ddrv19} solves independently the $N$ one-dimensional problems in \eqref{eq:model_decoupled}. Then the solution is visualized by stacking the one-dimensional reconstructions obtaining a two-dimensional image. Finally, the obtained results by the algorithm presented in this paper are depicted in the last row. 

\begin{figure}
\centering
\begin{minipage}{0.4\linewidth}
	\centering
	\includegraphics[width=\linewidth]{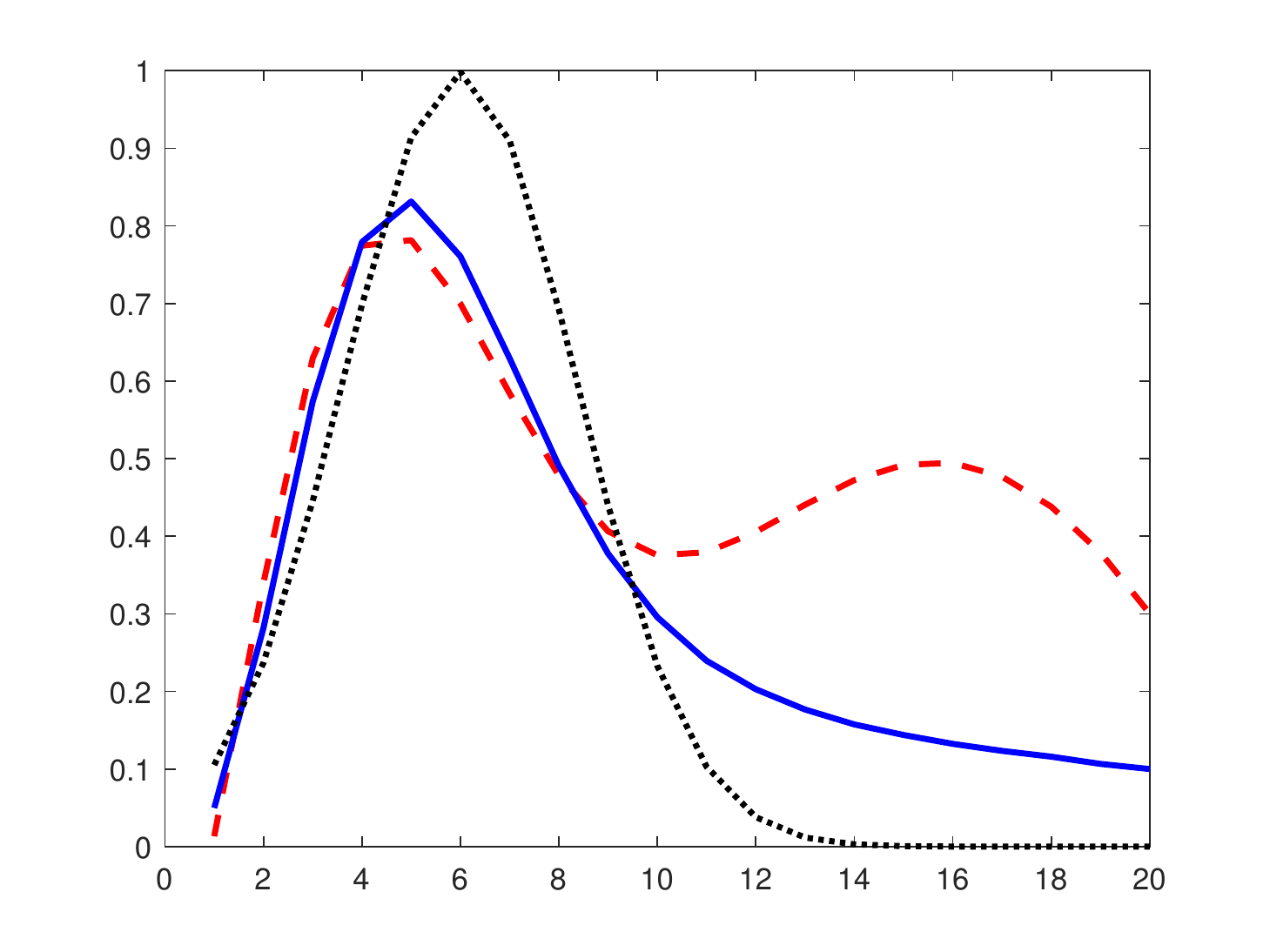}\\(a)
\end{minipage}
\begin{minipage}{0.4\linewidth}
	\centering
	\includegraphics[width=\linewidth]{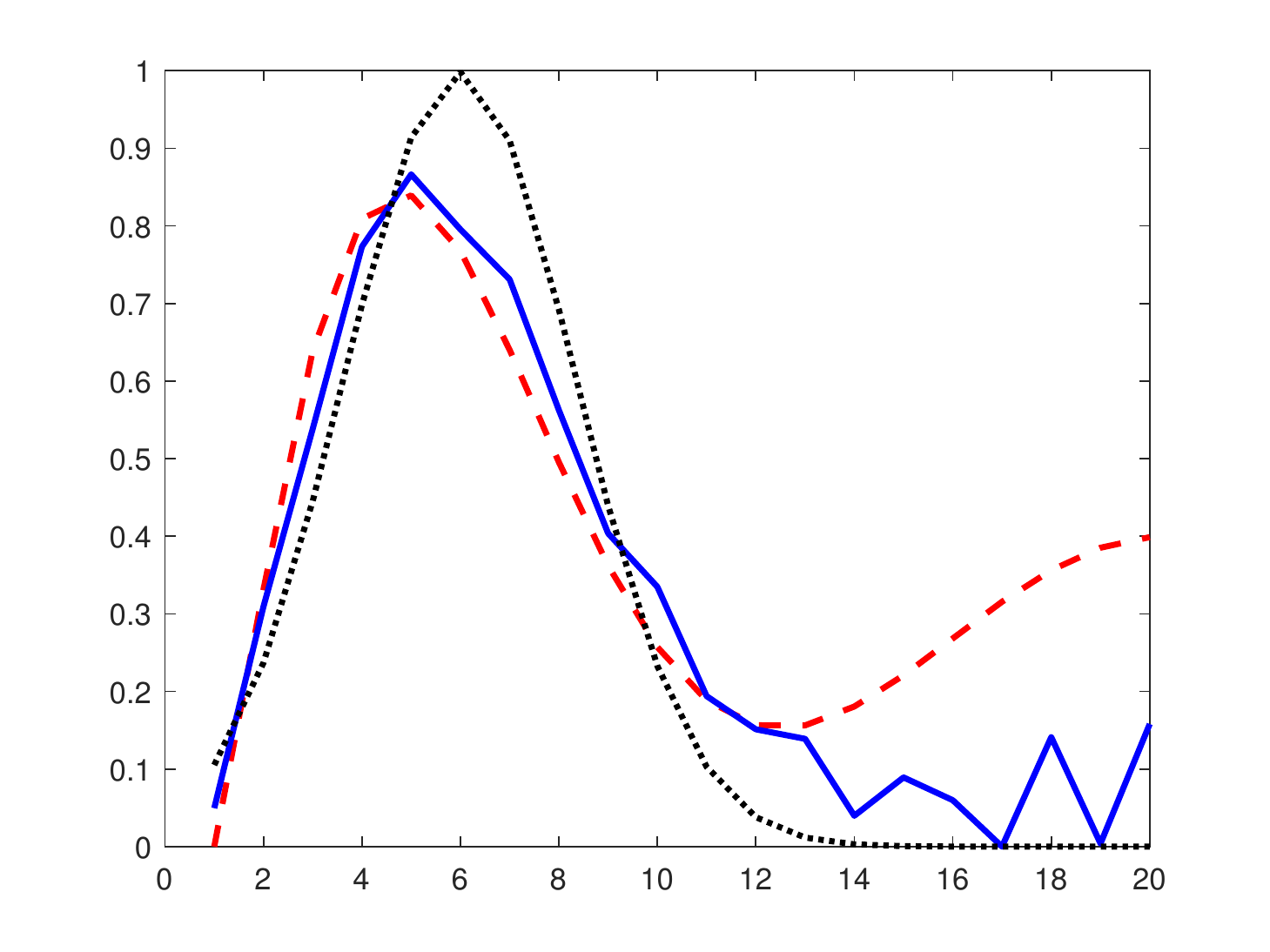}\\(b)
\end{minipage}
\caption{One-dimensional reconstruction of the two-dimensional synthetic model generated in Figure~\ref{Figure2}. On the left, $\boldsymbol{\sigma}_{16}$ in the model generated by Geophex - GEM 2 (panel (a)). On the right $\boldsymbol{\sigma}_{16}$ in the model generated by CDM Explorer (panel (b)). The exact solution is depicted by a black dotted line, the profile computed by the method in \cite{ddrv19} is reported with a red dashed line, and the profile obtained with our method is shown by a blue solid line.}
\label{Figure3}
\end{figure}

We can see that, for both the instruments configuration, the reconstruction of the electrical conductivity obtained from the variational model described here is much more accurate than the ones obtained by the method described in \cite{ddrv19}, in which there is not lateral continuity in the results. This shows that the additional regularization term in \eqref{eq:model} significantly helps the reconstruction of an approximate solution. In Figure~\ref{Figure3} we report the reconstruction of a single sounding (the
16th column of the two-dimensional synthetic model of Figure~\ref{Figure2} for both the devices). We compare the one-dimensional solution computed by the algorithm in \cite{ddrv19} to the one obtained with the proposed method. We can observe that the introduction of the regularization term allows us to improve the accuracy of the reconstruction. Finally, we compute the Relative Restoration Error (RRE) for each computed solution, defined by
$$
{\rm RRE}(\Sigmab)=\frac{\norm{\Sigmab-\Sigmab_{\rm exact}}_F}{\norm{\Sigmab_{\rm exact}}_F},
$$
where $\Sigmab_{\rm exact}$ denotes the exact solution of the problem. We report the obtained results in Table~\ref{tbl:rre}. We can observe that the RRE obtained with our proposal is significantly lower than the one obtained with the method proposed in \cite{ddrv19}.

\paragraph{Test 2.}
This example reports the results for different dimensions of the matrix solution: $50$, $100$ and $200$ soundings along a $10$ m straight line and $20$, $50$ and $100$ layers, respectively.
In this case, the synthetic data is generated by the Geophex GEM-2 ($\rho = 1.66~{\rm m},\; f = 775,\; 1175,\; 3925,\; 9825,\; 21725,\; 47025~{\rm Hz}$) with two orientations of the coils and one measurement height $h = 1$ m. We use the same parameters as in the previous example. 

\begin{figure}
\centering
\begin{minipage}{0.3\linewidth}
	\centering
	\includegraphics[width=\linewidth]{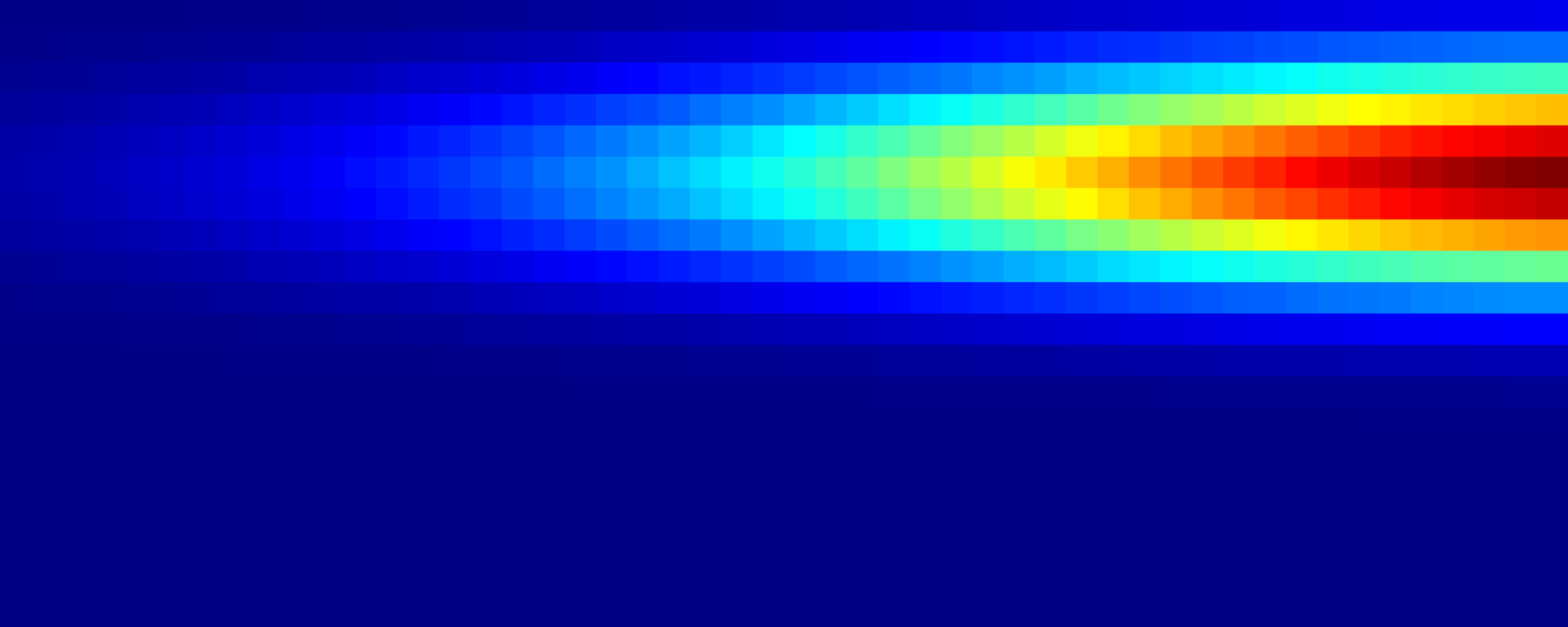}\\(a)
\end{minipage}
\begin{minipage}{0.3\linewidth}
	\centering
	\includegraphics[width=\linewidth]{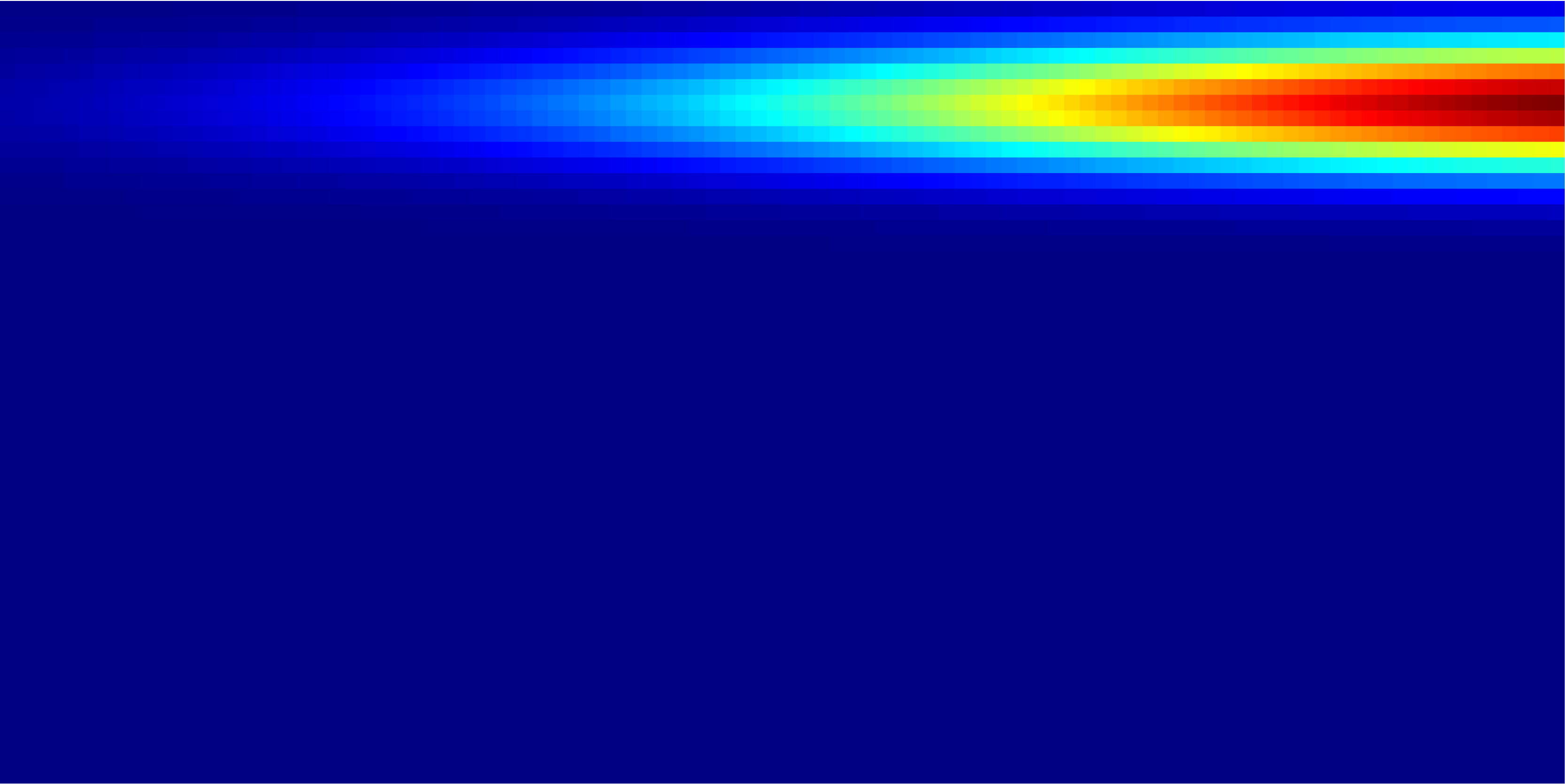}\\(b)
\end{minipage}
\begin{minipage}{0.3\linewidth}
	\centering
	\includegraphics[width=\linewidth]{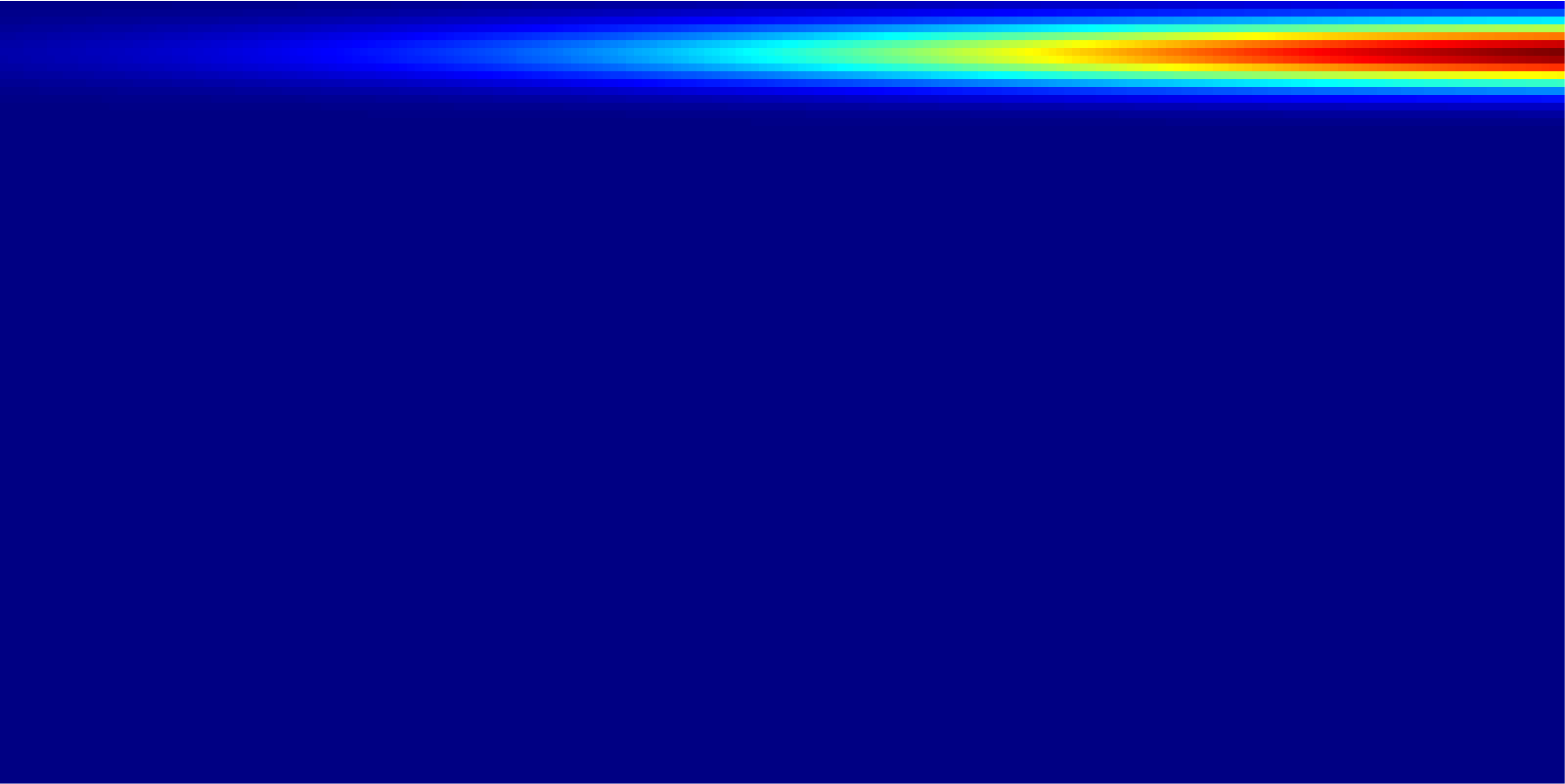}\\(c)
\end{minipage}

\begin{minipage}{0.3\linewidth}
	\centering
	\includegraphics[width=\linewidth]{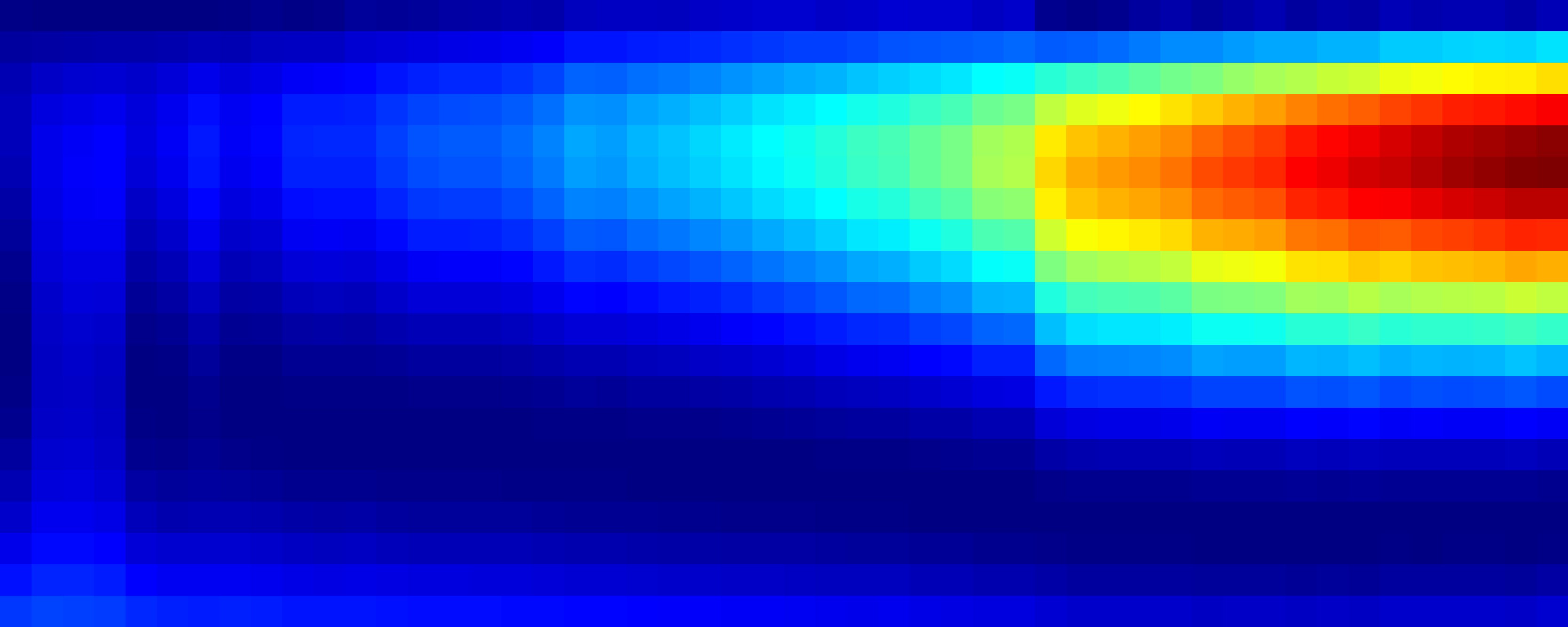}\\(d)
\end{minipage}
\begin{minipage}{0.3\linewidth}
	\centering
	\includegraphics[width=\linewidth]{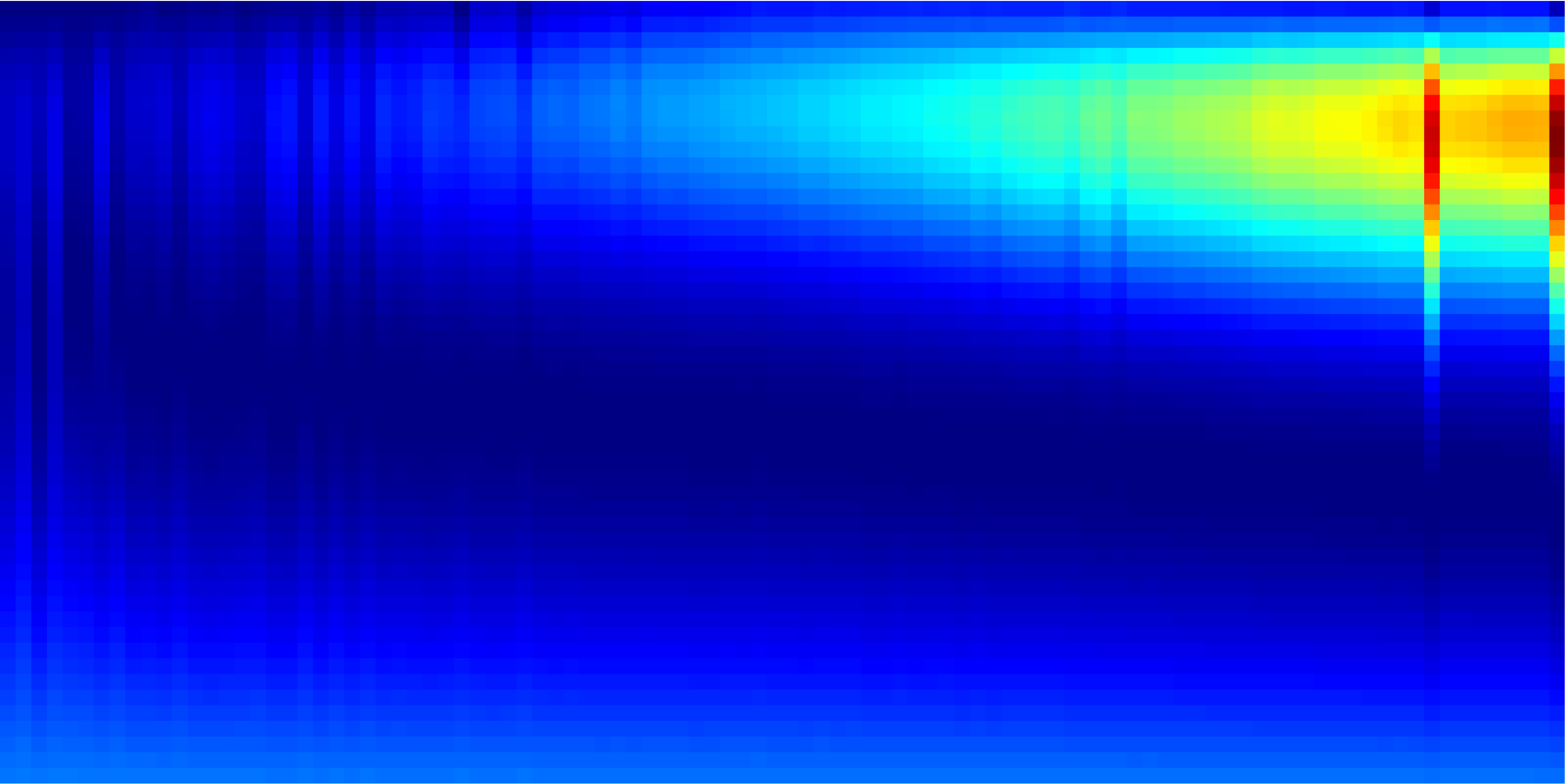}\\(e)
\end{minipage}
\begin{minipage}{0.3\linewidth}
	\centering
	\includegraphics[width=\linewidth]{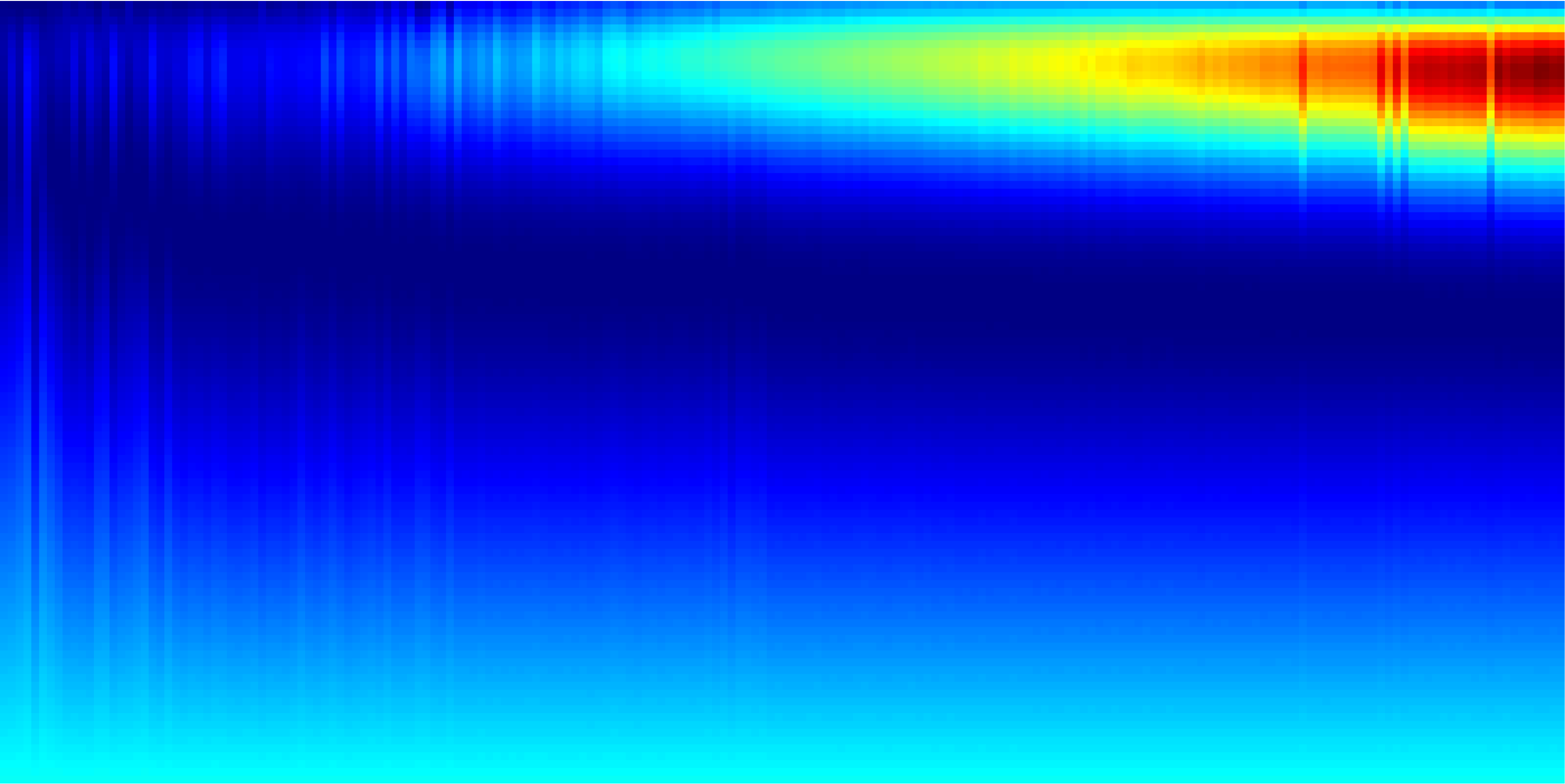}\\(f)
\end{minipage}

\begin{minipage}{0.3\linewidth}
	\centering
	\includegraphics[width=\linewidth]{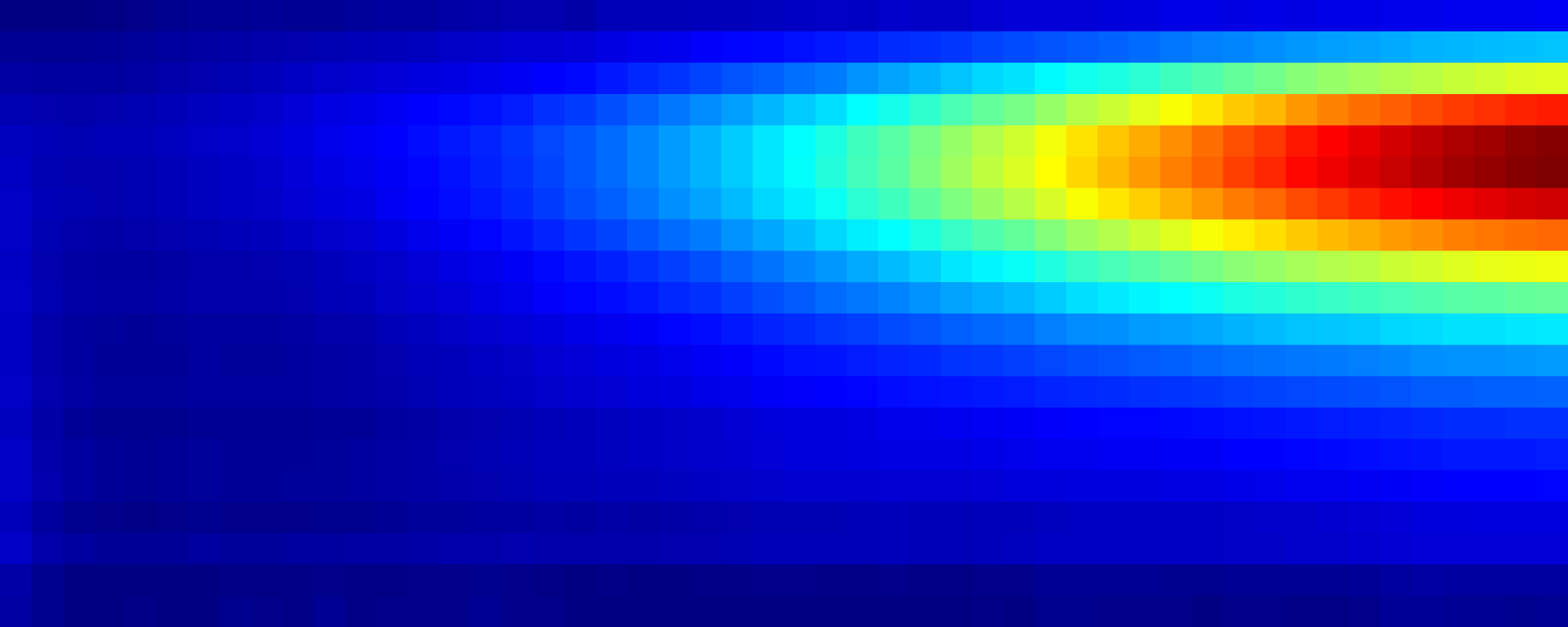}\\(g)
\end{minipage}
\begin{minipage}{0.3\linewidth}
	\centering
	\includegraphics[width=\linewidth]{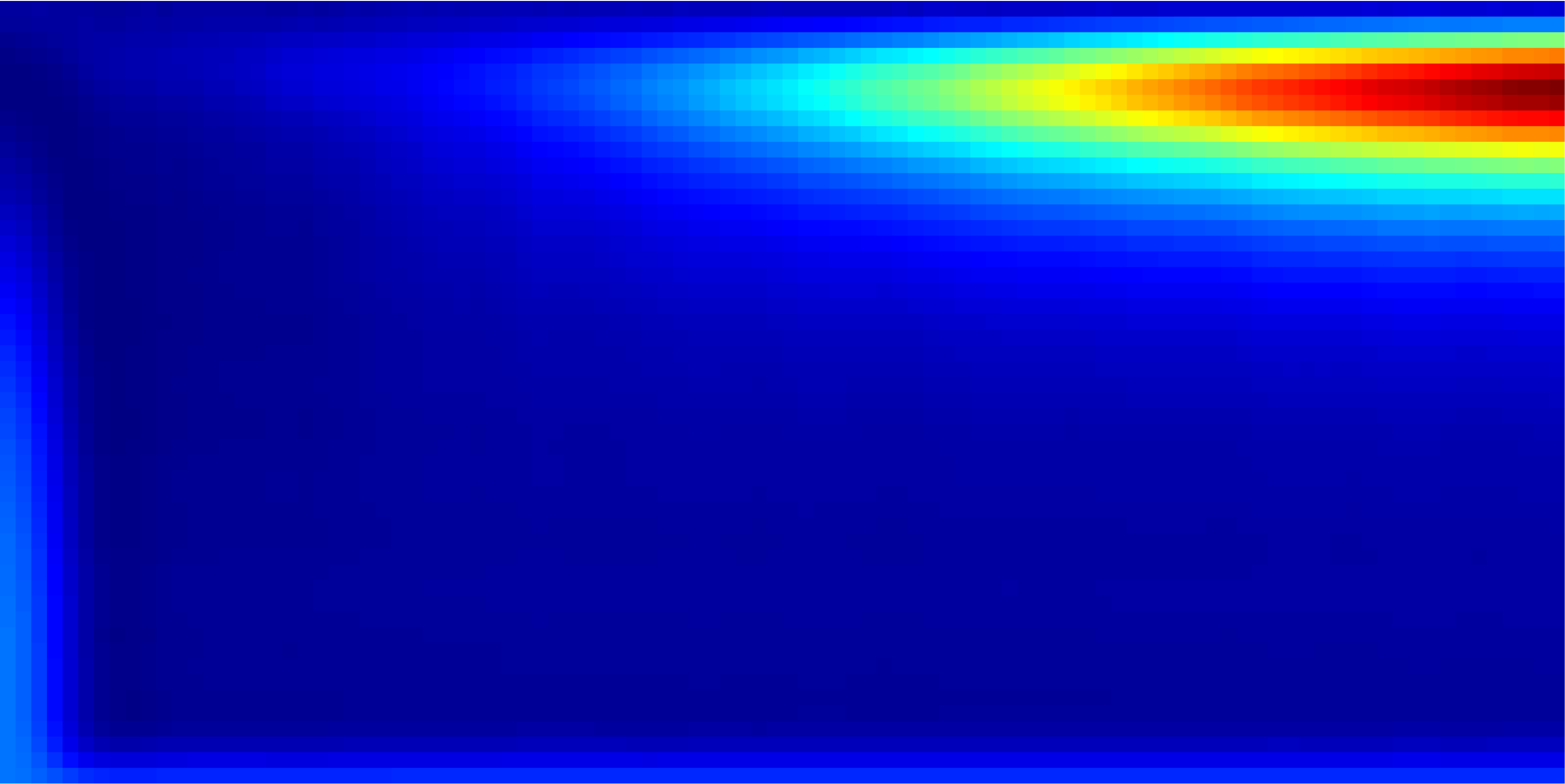}\\(h)
\end{minipage}
\begin{minipage}{0.3\linewidth}
	\centering
	\includegraphics[width=\linewidth]{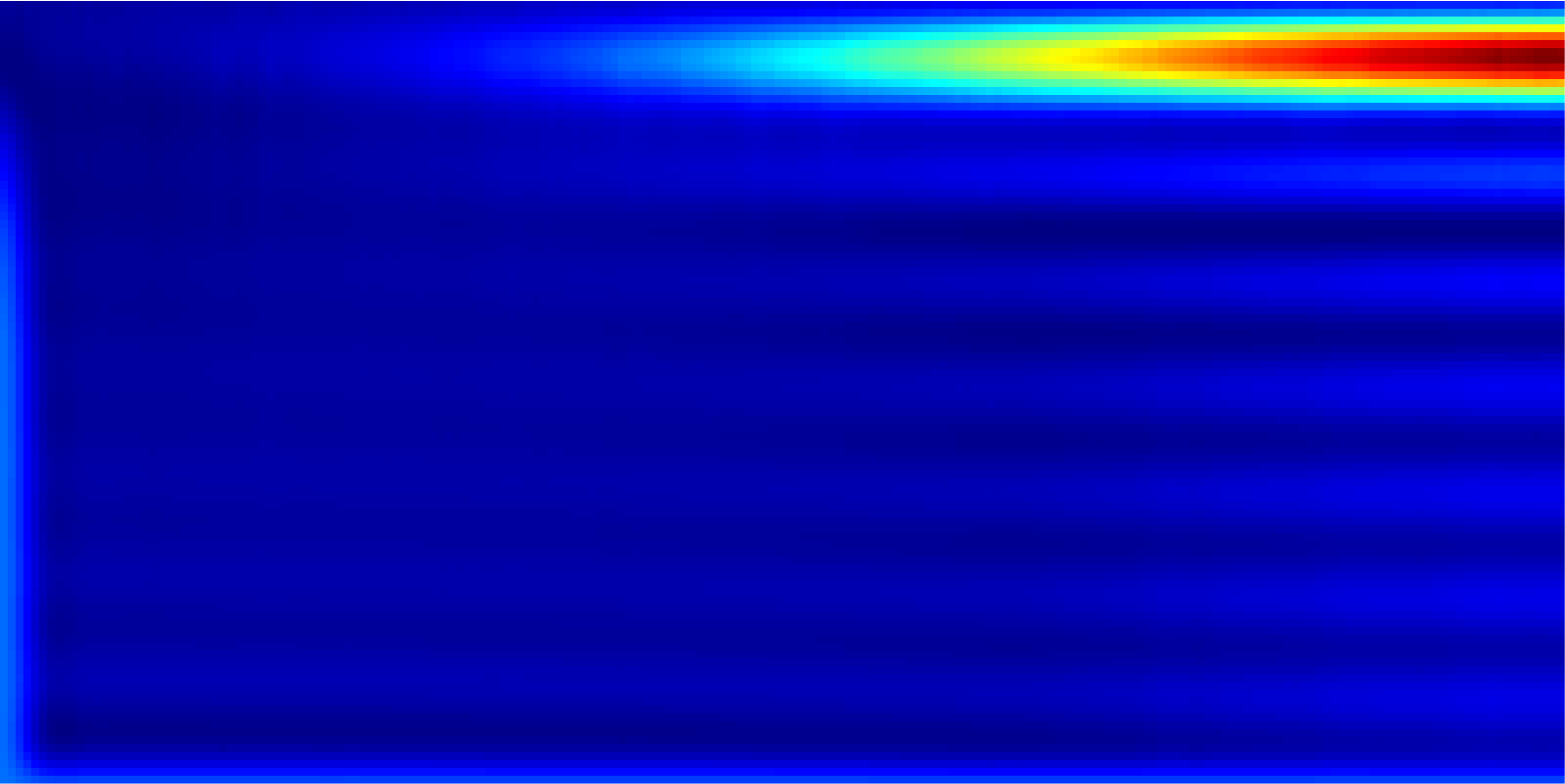}\\(i)
\end{minipage}

\caption{Reconstruction of the electrical conductivity from data generated by the Geophex - GEM 2. The first row reports the exact solutions, the second row shows the reconstructions obtained with the method in \cite{ddrv19}, and the third row collects the approximate solution computed by our algorithmic proposal. Each column reports the results obtained with different examples. In the first column we consider $\Sigmab_1\in\R^{20\times 50}$, in the second column $\Sigmab_2\in\R^{50\times 100}$, and in the third column $\Sigmab_3\in\R^{100\times 200}$. The data values are affected by Gaussian noise with noise level $\delta = 10^{-2}$, the maximum number of iterations is equal to $50$, $p=2$, $q=0.1$, and the regularization parameter $\gamma = 10^{-4}$.}
\label{Figure4}
\end{figure}

In the first column of Figure~\ref{Figure4}, we report the reconstructions of $\Sigmab_1\in\R^{20 \times 50}$. The second column depicts the reconstructions of $\Sigmab_2\in\R^{50 \times 100}$. The last column increases the dimension of the exact solution to $\Sigmab\in\R^{100 \times 200}$. In all the experiments the right-hand sides are affected by Gaussian noise with $\delta = 10^{-2}$, the maximum number of iterations is equal to $50$, $p=2$, $q=0.1$, the regularization parameter $\gamma = 10^{-4}$, and $\sigma_0=0.2$.

From the visual inspection of the reconstructions in Figure~\ref{Figure4} we can observe that the proposed model is able to capture more accurately the structure of the soil in all proposed examples. Moreover, we can see that, since the reconstruction is not obtained by simply stacking one-dimensional vectors, the approximate solutions computed by our proposal are smoother than the ones obtained with the method described in \cite{ddrv19}. This is confirmed by the values of the RRE reported in Table~\ref{tbl:rre}.

\begin{table}
	\caption{Relative Restoration Error (RRE) obtained in all the considered examples on synthetic data with the algorithm proposed in \cite{ddrv19} and our proposal.}
	\label{tbl:rre}
	\begin{center}
		\begin{tabular}{l|l|l|l}
			\multicolumn{2}{l}{Example}&Method&RRE\\\hline
			\multirow{4}{*}{Test 1}&\multirow{2}{*}{Geophex}&\cite{ddrv19}&$0.55439$\\
			&&Algorithm~\ref{algo:final}&$0.37832$\\\cline{2-4}
			&\multirow{2}{*}{Explorer}&\cite{ddrv19}&$0.41573$\\
			&&Algorithm~\ref{algo:final}&$0.35842$\\\hline
			\multirow{6}{*}{Test 2}&\multirow{2}{*}{$20\times50$}&\cite{ddrv19}&$0.28644$\\
			&&Algorithm~\ref{algo:final}&$0.25646$\\\cline{2-4}
			&\multirow{2}{*}{$50\times100$}&\cite{ddrv19}&$0.65784$\\
			&&Algorithm~\ref{algo:final}&$0.36258$\\\cline{2-4}
			&\multirow{2}{*}{$100\times200$}&\cite{ddrv19}&$0.93571$\\
			&&Algorithm~\ref{algo:final}&$0.36258$\\\hline
		\end{tabular}
	\end{center}
\end{table}
\subsection{Experimental data}
In this section, we consider an experimental data set collected with a multiconfiguration EMI device at the Molentargius Saline Regional Nature Park, located east of Cagliari in southern Sardinia, Italy. The dataset was first studied in \cite{ddrv19}; see \cite{ddrv19} for a description of the geographical location of the test site.

As described in \cite{ddrv19}, the CDM Explorer ($\rho = 1.48, 2.82, 4.49~{\rm m}$ and $f = 10~{\rm kHz}$) was used to collect the EMI data along a $200$ m straight-line and was carried at $0.9$ m of height above the ground, providing three simultaneous measurements for each orientation of the device. The data were collected in continuous mode, with a $0.5$ s time step, first using the horizontal orientation and then the vertical one. With the aim of having the same number of equally spaced measurement points, the authors of \cite{ddrv19} merged the data by spatially resampling at $0.5$ m intervals from a common starting point. This allowed to set up a dataset consisting of a series of $401$ depth soundings with six complex (quadrature and in-phase components) measurements each, suitable for recovering the soil electrical conductivities in order to reconstuct the water table. The data set is available at the web page \url{http://bugs.unica.it/cana/datasets/}.

\begin{figure}
\centering
\begin{minipage}{0.45\linewidth}
	\centering
	\includegraphics[width=\linewidth]{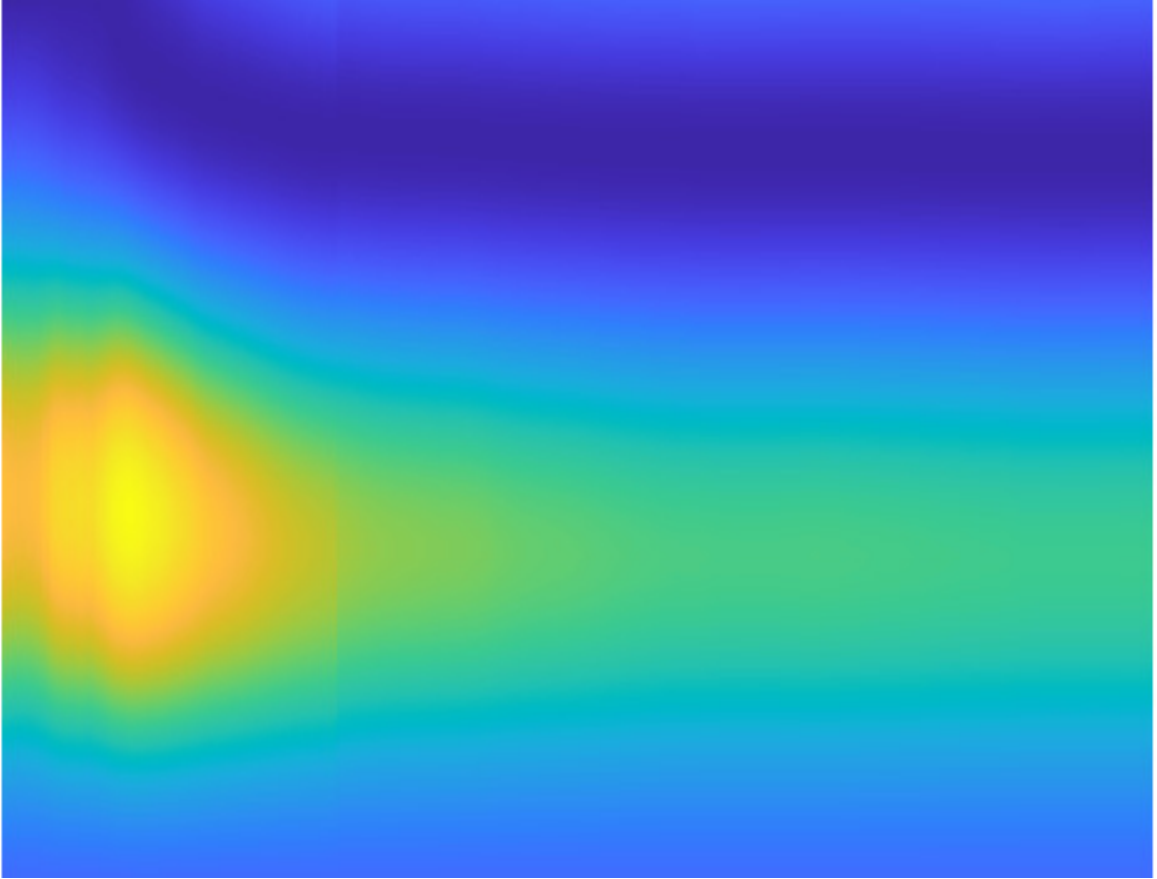}\\(a)
\end{minipage}
\begin{minipage}{0.45\linewidth}
	\centering
	\includegraphics[width=\linewidth]{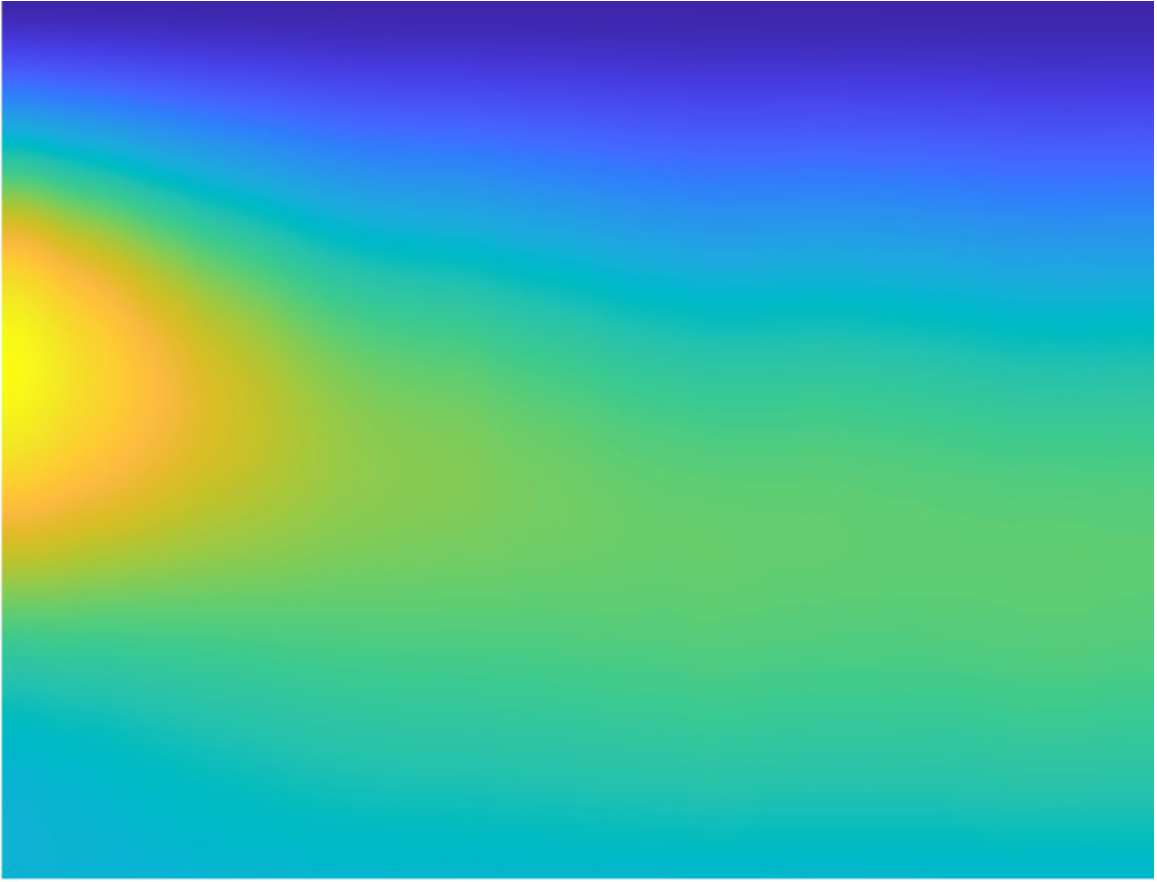}\\(b)
\end{minipage}
\caption{Reconstruction of the electrical conductivity from real data collected at Molentargius Saline Regional Nature Park in Cagliari by the CDM Explorer. The panel on the left reports the reconstruction obtained with the method in \cite{ddrv19}, and the right one collects the approximate solution computed by our algorithmic proposal. In both panels we consider $\Sigmab\in\R^{50\times 401}$ and a maximum number of iterations equal to $10^3$. In panel (a) we use the first derivative operator as the regularization matrix, while in panel (b) we fix $p=2$, $q=0.1$, and we set the regularization parameter $\gamma = 10^{-4}$.}
\label{figrealdata}
\end{figure}

In Figure~\ref{figrealdata}, we compare the results obtained by the method described in \cite{ddrv19} where we used the first derivative operator as the regularization matrix, with the ones computed by Algorithm~\ref{algo:final} setting $p=2$, $q=0.1$, and the regularization parameter $\gamma = 10^{-4}$. In both cases the maximum number of iterations has been fixed to $10^3$. We discretize the soil with $50$ layers leading to a reconstruction $\Sigmab\in\R^{50 \times 401}$. 

\begin{figure}
\centering
\begin{minipage}{0.45\linewidth}
	\centering
	\includegraphics[width=\linewidth]{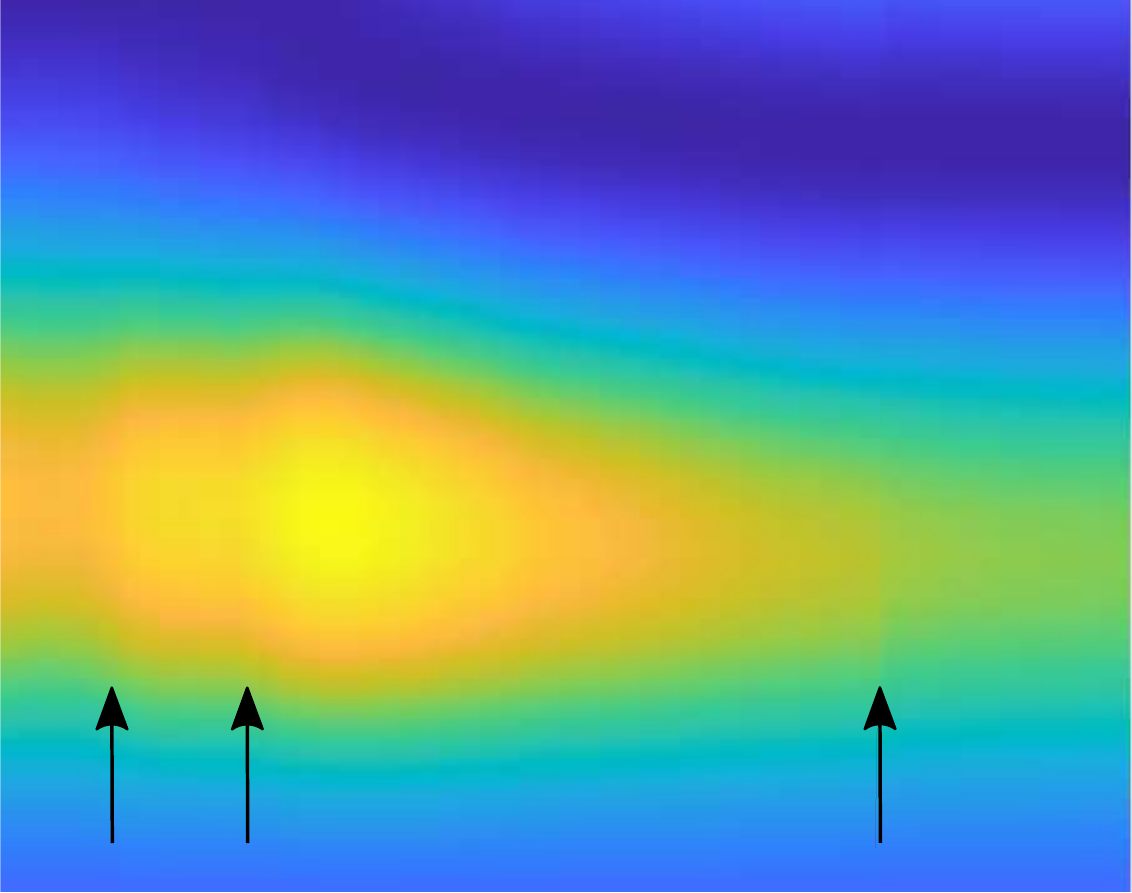}\\(a)
\end{minipage}
\begin{minipage}{0.45\linewidth}
	\centering
	\includegraphics[width=\linewidth]{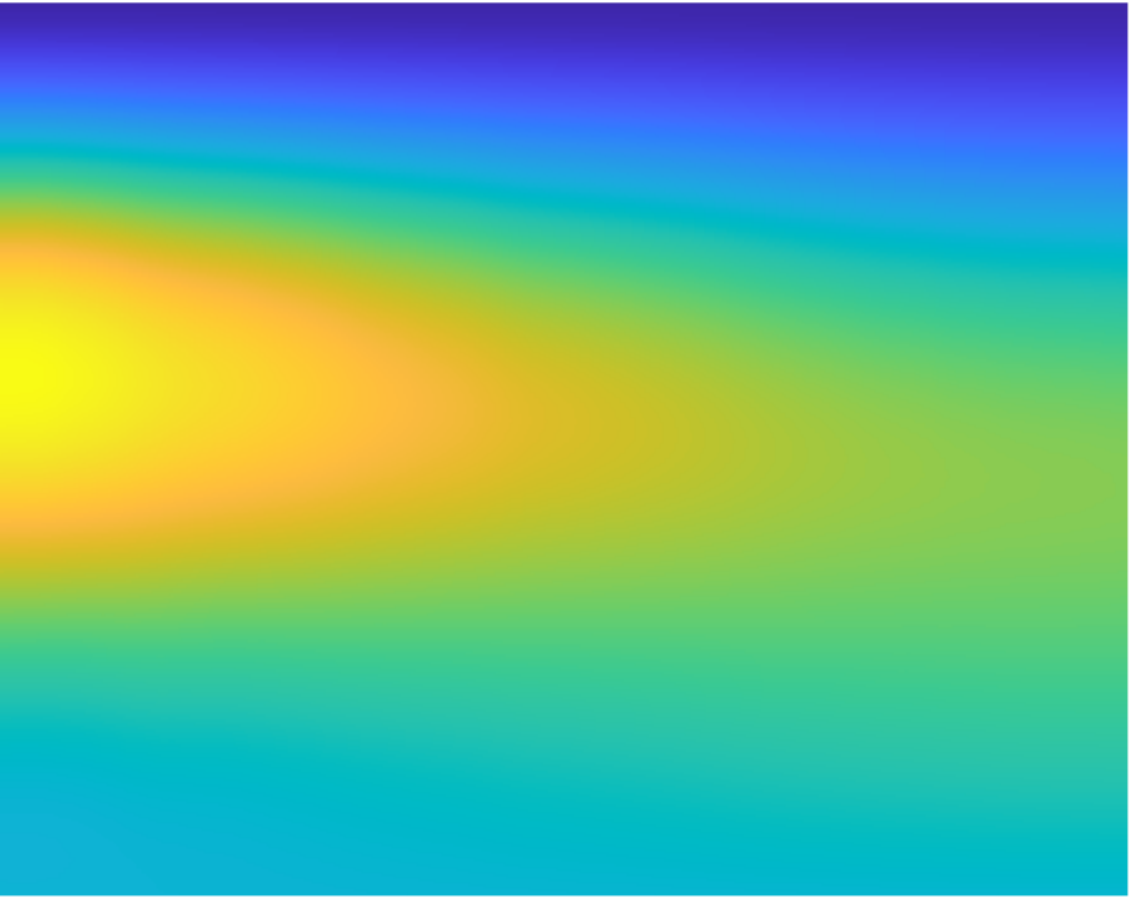}\\(b)
\end{minipage}
\caption{Plot of the first 100 soundings of the reconstructions in Figure~\ref{figrealdata}.}
\label{figrealdata2}
\end{figure}

From a graphic examination of the reconstructions in Figure~\ref{figrealdata}, we can observe that the proposed method (Figure~\ref{figrealdata}(b)) is able to capture more accurately the structure of the soil, avoiding any ``splicing''. On the other hand the reconstruction obtained with the algorithm in \cite{ddrv19}, shown in Figure~\ref{figrealdata}(a), is ``spliced'' in several points, that are highlighted in Figure~\ref{figrealdata2}. Moreover, we can observe that both reconstructions have the same shape and structure, showing the water table interface, even if it is not easy to identify its exact depth.

We remark here that in these experiments, the same starting model was adopted for both of the reconstructions but still, there is an important dependence on the initial solution of the iterative method, also in our optimization algorithm.

\section{Conclusions}\label{sect:concl}
In this paper we have proposed a non-linear model for the inversion of FDEM data. Although the two-dimensional problem can be seen as a sequence of  ``stacked'' independent one-dimensional problems, the proposed model couples them. This allows us to exploit the ``horizontal'' information to largely improve the quality of the computed solution avoiding the ``splicing'' that occurs when we deal with noisy data.

We were able to show that the outlined minimization problem has a solution, albeit non-unique, and that it induces a regularization method. We provided an algorithm for the solution of the problem and showed its performances on some synthetic and real data. From these experiments we were able to show that our method reliably computes accurate solutions even when compared with state-of-the-art method like the one in \cite{ddrv19}.

Matters of future research include an efficient implementation exploiting Krylov subspaces. Moreover, we plan on using the Alternating Direction Multiplier Method (ADMM) and its accelerations (see, e.g., \cite{boyd,GOSB14,BDD20}) for efficiently computing a solution of \eqref{eq:model}.

\section*{Acknowledgements}
The authors would like to thank Prof. G. Deidda for providing us with the experimental data and, more in general, for his help, Prof. G. Rodriguez for the precious discussions and  the anonymous referees for their comments that greatly improved the quality of this paper. The authors are members of the GNCS group of INdAM. Alessandro Buccini is partially supported by Regione Autonoma della Sardegna research project ``Algorithms and Models for Imaging Science [AMIS]'' (RASSR57257, intervento finanziato con risorse FSC 2014-2020 - Patto per lo Sviluppo della Regione Sardegna).
Patricia D\'iaz de Alba is partially supported by INdAM-GNCS 2020 project ``Tecniche numeriche per l'analisi delle reti complesse e lo studio dei problemi inversi''.

\bibliographystyle{siam}
\bibliography{bibliography}

\end{document}